\documentclass[a4paper,12pt]{article}
\usepackage[utf8]{inputenc}
\usepackage[T1]{fontenc}
\usepackage[english]{babel}
\usepackage{mathtools,amssymb,amsthm}
\usepackage{hyperref}
\usepackage{enumitem}
\usepackage{esint}
\usepackage{graphicx}
\usepackage{lmodern}
\usepackage{geometry} 
\usepackage{amsbsy} 
\usepackage{mathrsfs}
\usepackage{authblk}
\hypersetup{
  colorlinks=true,
  linkcolor=blue,      % internal links: sections, theorems, equations
  citecolor=blue,      % \cite
  urlcolor=blue,       % \url, \href
  filecolor=blue,
  breaklinks=true,     % allow links to break across lines
  bookmarksnumbered=true,
  pdftitle={Regularity of planar Griffith  minimizers with connected cracks through a monotonicity formula},
  pdfauthor={Labourie, Lemenant, Nahon},
}

\newcommand{\A}{\mathbb{A}}

\newcommand \Om{\Omega}
\newcommand \om{\omega}
\newcommand \eps{\varepsilon}
\newcommand \Cy{\mathcal{C}}

\newcommand \Ha{\mathcal{H}^1}
\newcommand \G{\mathcal{G}}
\newcommand \B{B}

\newcommand \N{\mathbb{N}}
\newcommand \Z{\mathbb{Z}}
\newcommand \C{\mathbb{C}}
\newcommand \R{\mathbb{R}}
\newcommand \loc{\mathrm{loc}}

\newcommand{\HH}{\mathcal{H}} % mesure de Hausdorff

\newcommand{\CC}[1]{\mathbb{A} #1:#1}
\usepackage{braket} % '\set', '\Set'
\renewcommand{\H}{\mathcal{H}}

\usepackage[breakable]{tcolorbox}
\usepackage{xcolor}
\usepackage{soul}
\definecolor{mycolor}{RGB}{229,228,200}
\definecolor{pinkcolor}{RGB}{247,202,193}

\newtcolorbox{changebox}{
    colback=mycolor,    % Couleur de fond (15% gris)
    colframe=mycolor,   % Couleur du cadre
    boxrule=0pt,
    sharp corners,      % Coins carrés
    left=0pt,
    right=0pt,
    top=0pt,
    bottom=0pt,
    breakable}

\numberwithin{equation}{section}

\newtheorem{theorem}{Theorem} 
\newtheorem{corollary}[theorem]{Corollary} 
\newtheorem{lemma} [theorem]{Lemma} 
\newtheorem{remark} [theorem]{Remark} 
\newtheorem{proposition} [theorem]{Proposition} 
 
\theoremstyle{definition} 
\newtheorem{definition}[theorem] {Definition}

\title{Regularity of planar Griffith  minimizers with connected cracks through a monotonicity formula}
\author[1]{Camille Labourie}
\author[1]{Antoine Lemenant}
\author[2]{Mickaël Nahon}

\affil[1]{\small Université de Lorraine, CNRS, IECL, F-54000 Nancy, France}
\affil[2]{\small Univ. Grenoble Alpes, CNRS, Grenoble INP, LJK, 38000 Grenoble, France}
\date{}

\begin{document}
\maketitle
\abstract{We prove a new monotonicity formula for the two-dimensional Neumann elasticity equation with  connected boundary. From this we deduce a   regularity result for two-dimensional  minimizers of the Griffith energy with  connected crack set: the crack may be partitioned into regular points, triple junctions, and cracktips. We give a classification of the global  minimizers of the Griffith energy with connected crack set. Finally, we discuss some consequences of the monotonicity formula on the evolution of the stress intensity factors along a crack, which may present independent interest.}

\tableofcontents
\newpage
\section{Introduction}

Given some Lamé parameters $\lambda,\mu$, such that
$$\mu>0,\ \lambda+\mu>0,$$
we define, for any open set $\Omega \subset\R^2$, for any relatively closed set $K\subset \Omega$, and for any vector field $u\in H^1_{\loc}(\Omega\setminus K,\R^2)$, the Griffith energy
$$\G(u,K|\Omega ):=\int_{\Omega \setminus K}\A e(u):e(u)+\Ha(\Omega\cap K),$$
where $\mathcal{H}^1$ is the one-dimensional Hausdorff measure, $e(u)=\frac{Du + (Du)^*}{2}$ and $\A$ is the Lamé tensor defined by
\begin{equation}\label{eq_deftensorA}
\A\xi=\lambda \mathrm{tr}(\xi)I_2+2\mu\xi.
\end{equation}

The Griffith energy provides one of the fundamental variational frameworks for describing crack propagation in brittle materials. In this approach, a configuration is characterized by a displacement $u$ and a crack set $K$, and its energy results from the competition between the elastic energy stored in the material and the cost associated with the creation of new fracture surfaces. In its modern variational formulation, the problem therefore consists in minimizing the energy $\G(u,K|\Omega )$ with a boundary datum $u=g$ on $\partial \Omega$.

The regularity theory for Griffith minimizers builds on the classical Mumford–Shah theory. In dimension two, Babadjian–Iurlano–Lemenant \cite{BIL} established partial regularity of the crack set, while Friedrich–Labourie–Stinson \cite{FLK} later proved an $\varepsilon$-regularity result yielding local $C^{1,\alpha}$ a.e. regularity and showed that the singular set has Hausdorff dimension strictly smaller than one. In arbitrary dimension, Labourie–Lemenant obtained $C^{1,\alpha}$ regularity under a flatness and separating condition, while   partial regularity in higher dimensions remains open.

In this paper, we establish the full regularity of every isolated connected component of the crack set \(K\). The main ingredient is a new monotonicity formula for the two-dimensional elasticity equation with Neumann condition, which enables us to characterize the blow-up limits of minimizers. In this respect, our work can be viewed as a counterpart of Bonnet's result \cite{b} on the minimizers of the Griffith energy. We recall that the existence of blow-up limits was already established in \cite{lemenant_labourie} using the notion of uniform concentration. More precisely, it was proved therein that, up to a subsequence, every blow-up sequence converges to a global minimizer.

In order to present the results of this paper, let us give some more precise definitions.

\begin{definition}[Admissible pair]
An \emph{admissible pair} $(u,K)$ consists of a relatively closed subset $K\subset\Omega$ and a vector
            field $u \in H^1_{\mathrm{loc}}(\Omega \setminus K,\R^2)$ such that 
            \begin{equation}
                \H^1(K \cap \B_{x,r})>0, \quad\forall x \in K,\,\forall r>0.
            \end{equation}
and $\G(u,K|V)<+\infty$ for any compactly contained open set $V\Subset\Om$. Here $B_{x,r}$ is the open ball of $\R^2$ of center $x$, of radius $r$.
 \end{definition}

 \begin{definition}[Topological Competitor]\label{def_topcompetitor}
 Let $(u,K)$ be an admissible pair. Let $x_0 \in \Omega$ and $r > 0$ be such that $\overline{\B_{x_0,r}} \subset \Omega$.
            A \emph{topological competitor} of $(u,K)$ in $\B_{x_0,r}$ is an admissible pair $(v,L)$ such that
            \begin{equation}
                L \setminus \B_{x_0,r} = K \setminus \B_{x_0,r} \quad \text{and} \quad v = u \quad \text{a.e. in } \Omega \setminus \left(K \cup \B_{x_0,r}\right).
            \end{equation}
            and
            \begin{gather}\label{eq_topo_competitor}
                \text{all points $x,y \in \Omega \setminus (K \cup \B_{x_0,r})$ that are separated by $K$}\\
                \text{are also separated by $L$}.
            \end{gather}
            In other words, if $x,y \in \Omega \setminus (K \cup \B_{x_0,r})$ belong to different connected components of $\Omega \setminus K$, they also belong to different connected components of $\Omega \setminus L$.
\end{definition}

 Now here is the notion of minimizers.
 
\begin{definition}[Topological minimizers]\label{defAlmostMin} A \emph{topological minimizer} of the Griffith functional in $\Om$ is an admissible  pair $(u,K)$ such that for all $x \in \Omega$ and $r > 0$ with $\overline{\B_{x,r}}
    \subset \Omega$ and for any topological competitor $(v,L)$ of $(u,K)$ in $\B_{x,r}$, we have
$$ \G(u,K| \B_{x,r} )\leq \G(v,L| \B_{x,r} ).$$
\end{definition} 

We will also consider the notion of \emph{global minimizer}, which arises as blow-up limits of almost minimizers as was proved in \cite{lemenant_labourie} . Here is the definition.

\begin{definition}[Global minimizers]
    A \emph{global minimizer}  $(u,K)$ is a topological minimizer with $\Omega=\R^2$.
\end{definition} 

Notice that \emph{topological} minimizers form a larger class than minimizers, since the class of admissible competitors is more restrictive. Historically, the topological constraint naturally arises in the blow-up procedure, as first observed by A. Bonnet \cite{b}. The main motivation for introducing topological minimizers is therefore to provide a unified framework encompassing both global minimizers arising as blow-up limits and local minimizers.

\medskip

Let us give some examples of possible global minimizers : the first two are known to be global topological minimizer, whereas the status of the third example is not established. We remind that a rigid motion is an affine function of the form $x\mapsto Ax+b$, for some antisymmetric matrix $A$.
\begin{itemize}
\item[i)](Line)  $K=\R\times\{0\}$, and there exist rigid motions $a,b:\R^2\to\R^2$ such that
\begin{equation}\label{eq_def_flat}
u=a\chi_{\R\times \R_+}+b\chi_{\R\times \R_-}.
\end{equation}
\item[ii)](Triple junction) $K=\left\{re^{i\theta}\text{ s.t. } r\geq 0,\ \theta\in\{0,\frac{2\pi}{3},\frac{4\pi}{3}\}\right\}$ and there exist rigid motions $a,b,c:\R^2\to\R^2$ such that
\begin{equation}\label{eq_def_triple}
u=a\chi_{\{0<\theta<\frac{2\pi}{3}\}}+b\chi_{\{\frac{2\pi}{3}<\theta<\frac{4\pi}{3}\}}+c\chi_{\{\frac{4\pi}{3}<\theta<2\pi\}}.
\end{equation}
\item[iii)](Cracktip) $K=]-\infty,0]\times\{0\}$ and there exists a rigid motion $a:\R^2\to\R^2$, such that for any $r>0,\theta\in (-\pi,\pi)$:
\begin{equation}\label{eq_def_cracktip}
u(re^{i\theta}) = a(re^{i\theta}) 
    \pm\sqrt{\frac{r}{16 \mu(\lambda+\mu)(\lambda+2\mu) \pi}}
    \begin{pmatrix}
        -(\lambda + \mu) \cos(3\theta /2) +(\lambda +5\mu)\cos(\theta/2)\\
        -( \lambda + \mu)\sin(3\theta /2) + (3\lambda +7\mu)\sin(\theta/2)
\end{pmatrix}.
\end{equation}
\end{itemize}

\subsection{Main results}

Our first main result is the classification of global minimizers.

\begin{theorem}\label{main1}
Let $(u,K)$ be a global  minimizer such that $K$ is a non-empty connected set. Then up to a translation and rotation, $(u,K)$ is of the form given in equation \eqref{eq_def_flat}, \eqref{eq_def_triple} or \eqref{eq_def_cracktip}.
\end{theorem}

A key ingredient in the proof of the theorem is a new monotonicity formula for the biharmonic energy, given in Theorem \ref{th_monotonie_disk}.

When \(K\) is a half-line, it is well known that the Lamé system admits two distinct $1/2$-homogeneous solutions, usually referred to as the \emph{opening} and \emph{sliding} modes. However, among the blow-up limits appearing in our classification, only the opening mode \eqref{eq_def_cracktip} occurs. Thus, our main result implies, in particular, that the sliding mode cannot arise as a blow-up limit. This exclusion is obtained through a careful examination of the elastic energy release along a slightly angled path (see Section~\ref{sec_stab}).

The classification of blow-up limits is then used to establish the full regularity of the crack, for an isolated connected component of a local minimizer. Here is our second main result.

\begin{theorem}\label{main2}
Let  $(u,K)$ be a topological  minimizer of the Griffith energy in $\Omega \subset \R^2$ such that $K$ is connected. Then there is a partition
$$K=\mathcal{R}\sqcup \mathcal{Y}\sqcup \mathcal{C}$$
where
\begin{itemize}
\item[(i) ] $\mathcal{R}$ is a locally finite union of $\mathcal{C}^1$ curves.
\item[(ii) ] $\mathcal{Y}$ is a discrete, relatively closed subset of $\Omega$, such that any blow-up of $(u,K)$ at a point of $\mathcal{Y}$ is a rotation of \eqref{eq_def_triple}. Moreover, for any $z\in\mathcal{Y}$, there exists $r>0$ such that $\B_{z,r}\cap K\setminus\{z\}\subset \mathcal{R}$ and there exists an homeomorphism $\Phi:\B_{z,r}\to \B_{r}$ such that $\Phi(\partial \B_{z,\rho})=\partial\B_{\rho}$ for every $\rho\in (0,r)$ and
$$\Phi(\B_{z,r}\cap K)=\left\{\rho e^{i\theta},\rho\in [0,r),\ \theta\in\{0,\frac{2\pi}{3},\frac{4\pi}{3}\}\right\}.$$
\item[(iii) ] $\mathcal{C}$ is a discrete, relatively closed subset of $\Omega$, such that any blow-up of $(u,K)$ at a point of $\mathcal{C}$ is a rotation of \eqref{eq_def_cracktip}. Moreover, for any $z\in\mathcal{C}$, there exists $r>0$ such that $\B_{z,r}\cap K\setminus\{z\}\subset \mathcal{R}$ and there exists an homeomorphism $\Phi:\B_{z,r}\to \B_{r}$ such that $\Phi(\partial \B_{z,\rho})=\partial\B_{\rho}$ for every $\rho\in (0,r)$ and
$$\Phi(\B_{z,r}\cap K)=(-r,0]\times\{0\}.$$
\end{itemize}
\end{theorem}

\subsection{Organization of the paper}
In section \ref{preliminaries} we give some further definitions and standard notions such as the Airy function and blow-up limits. 

In section \ref{sec_mono}, we develop a new monotonicity formula for the elasticity equation, that is stated in Theorem \ref{th_monotonie_disk} in terms of the Airy function. We also deduce from the monotonicity an explicit polynomial rate of decrease of the energy in balls of small radius, see Proposition \ref{prop_poly_growth}.

In section \ref{sec_stab}, we study the homogeneous minimizers of the Griffith problem. Using known results on the energy release rate of a crack with a small angle, we prove in Proposition \ref{prop_instability_slide} that a cracktip minimizer can only have a pure opening mode. We also obtain several new results on the existence of blow-up for a general solution of the elasticity equations (in Theorem \ref{th_existence_SIF}) and on the evolution of the stress intensity factors for a small perturbation (Theorem \ref{th_comparison_SIF} and Corollary \ref{cor_direction_kink}).

In section \ref{sec_epsreg}, we prove an $\varepsilon$-regularity result near flat solutions: any sufficiently flat topological minimizer is a smooth deformation of a straight fracture. The proof is based on the explicit polynomial decrease rate of energy proved in section \ref{sec_mono}.

Finally, in section \ref{sec_proof_main}, we prove the main results \ref{main1} and \ref{main2}, using both the $\varepsilon$-regularity result and the monotonicity formula.

\section{Preliminaries} \label{preliminaries}
\subsection{Notions on elasticity and fractures}
\noindent {\bf Notations.}
In all that follows, we let $\B_r$ be the centered ball of radius $r$ in $\R^2$, and $\B_{z,r}=z+\B_r$ for $z\in\R^2$. We will often assimilate $\R^2$ with $\C$.
We write $M_2(\R)$ the set of $2\times 2$ matrix, and for any $M\in M_2(\R)$ we write $\mathrm{Co}(M)$ the cofactor matrix of $M$. 
For any open set $\Om\subset\R^2$, for any relatively closed set $K\subset\Om$, for any $n\in\N^*$, we define 
$$H^{n}_{0,K}(\Om,\R)=\left\{\phi\in H^{n}(\Om,\R):\forall \eta\in\mathcal{C}^\infty_c(\Om),\ \eta\phi\in H^{n}_0(\Om\setminus K,\R)\right\}.$$
where $H^n_0(\Om\setminus K,\R)$ is the closure of $\mathcal{C}^\infty_c(\Om\setminus K,\R)$ for the $H^n$ norm. We also let
$$LD(\Om)=\left\{ u\in H^1_\loc(\Om,\R^2):e(u)\in L^2(\Om,S_2(\R))\right\}$$
where $S_2(\R)$ is the set of $2\times 2$ real symmetric matrices. We remind that \cite{chambolle} gives several density results of $H^1(\Om,\R^2)$ in $LD(\Om)$, but these two spaces are generally not equal.
\bigbreak
When $a$ and $b$ are two real quantities such that $b\geq 0$ and there exists a universal constant $C>0$ such that $a\leq Cb$, we write
$$a\lesssim b.$$
\noindent {\bf Convergence of pairs.}
The convergence of a sequence of admissible pairs $(u_i,K_i)$ to a limit $(u,K)$, is defined as follows
\begin{definition}\label{definition_conv}
Let $(u_i,K_i)_{i\geq 1}$ be a sequence of admissible pairs in a sequence of open sets $(\Omega_i)_{i\geq 1}$, and $(u,K)$ an admissible pair in an open set $\Omega$. Assume that for all compact set $H\subset \Omega$, we have $H\subset \Omega_i$ for all large enough $i$. We say that $(u_i,K_i)_{i\geq 1}$ converges to the pair $(u,K)$ in $\Omega$ when the following conditions are verified
\begin{enumerate}[label=(\roman*)]
\item $(K_i)_{i\geq 1}$ converges to $K$ in local Hausdorff distance: for all compact subset $H\Subset \Omega$ 
\begin{equation*}
\lim_{i \to +\infty} \left(\sup_{z \in K_i \cap H} \mathrm{dist}(z,K) + \sup_{z \in K \cap H} \mathrm{dist}(z,K_i)\right) = 0.
\end{equation*}
\item for all connected component $\mathcal{O}$ of $\Omega \setminus K$, there exists a sequence of rigid motions $(a_i^{\mathcal{O}})_{i\geq 1}$ such that for all compact subset $H \subset \mathcal{O}$,
\begin{equation*}
\lim_{i \to +\infty} \int_H |u_i - a_i^{\mathcal{O}} - u|^2 \; dx = 0.
\end{equation*}
\end{enumerate}
\end{definition}

This is the vectorial analogue of the convergence considered by {Bonnet} \cite{b}.
Note that the limit displacement $u$ is only determined up to a rigid motion in each connected component of $\Omega \setminus K$.

Topological minimizers are uniformly bounded in balls that meet the fracture set $K$, and this is a consequence of the following standard estimate.
\begin{proposition}\label{prop_Ahlfors}
There exists a constant $C_A>0$ such that the following holds. For any topological minimizer $(u,K)$ in $\Omega $, and for any ball $\B_{x,r}\Subset \Omega$ such that $\B_{x,r}\cap K\neq\emptyset$, then
$$\int_{\B_{x,r}}|e(u)|^2+\mathcal{H}^{1}(K\cap B_{x,r})\leq C_A r.$$
\end{proposition}
Indeed this is obtained by considering the competitor $(u1_{\B_{x,r}^c},K\cup \partial\B_{x,r}\setminus \B_{x,r})$.

As a consequence any sequence of minimizers $(u_i,K_i)$ admits a converging subsequence in the sense of definition \ref{definition_conv} ; we report here the result of \cite[Lemma 3.3. and Theorem 2.7]{lemenant_labourie}.

\begin{proposition}\label{prop_compactness}
Let $(u_i,K_i)_{i\geq 1}$ be a sequence of topological minimizers in a sequence of open set $(\Omega_i)_{i\geq 1}$. Let $\Omega$ be an open set such that for all compact set $H\subset \Omega$, we have $H\subset \Omega_i$ for all large enough $i$.

Then there exists a subsequence $(u_i,K_i)_{i\in I}$ (for some infinite set $I\subset\N^*$) and a minimizer $(u_\infty,K_\infty)$ in $\Omega$, such that $(u_i,K_i)_{i\in I}$ converges to $(u_\infty,K_\infty)$ in the sense of definition \ref{definition_conv}. Moreover, for all open subset $\om \Subset \Omega$, we have
    \begin{equation}\label{eq_blowup_energy}
        \int_{\om} \CC{e(u_{\infty})} \; dx = \lim_{i(\in I) \to +\infty} \int_{\om} \CC{e(u_i)}  \; dx
    \end{equation}
    and
    \begin{align*}
        \HH^{1}(K_{\infty} \cap \om) &\leq \liminf_{i(\in I) \to +\infty}  \HH^{1}(K_i \cap \om)\\
        \HH^{1}(K_{\infty} \cap \overline{\om}) &\geq \limsup_{i(\in I) \to +\infty}  \HH^{1}(K_i \cap \om).
    \end{align*}
    If in addition $\lim_{i(\in I)\to +\infty} \int_{\om}\A e(u_i):e(u_i) = 0$, then $u_\infty$ is piecewise rigid and $K$ is a minimal set in $\om$, under the topological constraint of definition \ref{def_topcompetitor}.
\end{proposition}

\begin{definition}[Standard rescaling of minimizers]\label{def_scaling}
If $(u,K)$ is a minimizer   in a domain $\Omega$, then for any $x\in\R^2$, $r>0$, the pair 
    
$$\left(u_{x,r},K_{x,r}\right)=\left(\frac{u(x+r\cdot)}{r^{1/2}},\frac{K-x}{r}\right)$$
is a minimizer in $\frac{\Omega-x}{r}$.
When $x=0$, we simply write $(u_r,K_r)$.
\end{definition}

\noindent{\bf Flatness.}
Let $(u,K)$ be a pair in an open set $\Omega$.
For any $x_0 \in K$ and $r_0 > 0$ such that $\B_{x_0,r_0} \subset \Omega$, we define the \textbf{bilateral} flatness $\beta_K(x_0,r_0)$ of $K$ in $\B_{x_0,r_0}$ as
\begin{equation*}
    \beta_K(x_0,r_0) := \frac{1}{r_0}\inf_{L} \max\left(\sup_{x \in K \cap \B_{x_0,r_0}} \mathrm{dist}(x,L),\sup_{x \in L \cap \B_{x_0,r_0}} \mathrm{dist}(x,K)\right),
\end{equation*}
where $L$ runs among affine lines passing through $x_0$. The flatness is invariant under rescaling, in the sense that
$$\beta_{K_{x_0,r_0}}(0,1)=\beta_{\frac{K-x_0}{r_0}}(0,1)=\beta_K(x_0,r_0)$$

There always exists an affine line $L$ which achieves the infimum. When there is no ambiguity on the considered set $K$, we write $\beta$ instead of $\beta_K$.

\medskip

\noindent{\bf Normalized elastic energy.}
For any $x_0 \in \Omega$ and $r_0 > 0$ such that $\B_{x_0,r_0} \subset \Omega$, we define the \emph{normalized elastic energy} of $u$ in $\B_{x_0,r_0}$ as
\begin{equation*}
    \omega_{u,K}(x_0,r_0) := \frac{1}{r_0} \int_{\B_{x_0,r_0} \setminus K} |e(u)|^2 \;dx.
\end{equation*}
As previously, we have the invariance
$$\om_{\frac{u(x_0+r_0\cdot)}{r_0^{1/2}},\frac{K-x_0}{r_0}}(0,1)=\om_{u,K}(x_0,r_0) $$

\subsection{Definition of the Airy function}

A central object from which we deduce a monotonicity formula for the elasticity equation, is the Airy function. We recall the notation 
$$\sigma:=\A e(u)=\lambda ({\rm tr}e(u)) I_2+2 \mu e(u). $$
 The purpose of this section is to prove the following proposition. 

\begin{proposition}\label{prop_def_airy}(Definition of the Airy function)
    Let $(u,K)$ be a global  minimizer such that $K$ is connected. Then there exists a scalar function $\phi : \R^2\to \R$ such that 
    \begin{enumerate}
        \item  $\phi\in H^2_{0,K}(\B_r)$ for all $r>0$,
        \item $\Delta^2 \phi = 0$ in $\R^2\setminus K$
        \item pointwisely in $\R^2\setminus K$ it holds, 
        $$D^2\phi=\mathrm{Co}(\sigma).
        $$
    \end{enumerate}
\end{proposition}

\begin{remark}
This proposition in fact only uses the minimality with respect to the function $u$ i.e. $u$ is a solution of the elasticity equation with homogeneous Neumann conditions along the crack.
\end{remark}

\begin{proof}We assume that $K$ is non-empty and more particularly that the origin belongs to $K$. 

The construction of the Airy function is rather classical. We can follow for instance the approach given in \cite[Proposition 4.3]{BCL1} (which itself relies on a construction already given by Chambolle in \cite{chambolle}). Our only concern is that the construction in \cite{BCL1} is local whereas here we want a function defined globally on $\R^2$.

For that purpose we consider the sequence of balls $\B_k$ for $k \in \N$ and we will construct $\phi$ inductively on $k$. We start applying \cite[Proposition 4.3]{BCL1} in the first ball $\B_1$ which yields a first function $\phi^{(1)}$ satisfying:
  \begin{enumerate}
        \item  $\phi^{(1)}\in H^2_{0,K}(\B_1)$,
        \item $\Delta^2 \phi^{(1)} = 0$ in $\B_1\setminus K$
        \item  $D^2\phi^{(1)}=\mathrm{Co}(\sigma)$ in $\B_1\setminus K$.
    \end{enumerate}
Then if $\phi^{(k)}$ is defined on $\B_k$ we show that one can extend it on $\B_{k+1}$ as follows: by application of \cite[Proposition 4.3]{BCL1} we get a function $\phi^{(k+1)}$ that satisfies:
  \begin{enumerate}
        \item  $\phi^{(k+1)}\in H^2_{0,K}(\B_{k+1})$,
        \item $\Delta^2 \phi^{(k+1)} = 0$ in $\B_{k+1}\setminus K$
       \item  $D^2\phi^{(k+1)}=\mathrm{Co}(\sigma)$ in $\B_{k+1}\setminus K$.
    \end{enumerate}
In particular we have $D^2(\phi^{(k+1)}-\phi^{(k)})=0$ in $ \B_{k}\setminus K$ : we deduce that on every connected component $\mathcal{O}$ of $\B_{k}\setminus K$ there exists an affine function $a^\mathcal{O}$ such that $\phi^{(k+1)}-\phi^{(k)}=a^\mathcal{O}$ in $\mathcal{O}$. Since both $\phi^{(k)}$ and $\phi^{(k+1)}$ belong to $H^2_{0,K}(\B_{k})$, we know that $a^\mathcal{O}=0$ $\mathrm{Cap}_{2,2}$-a.e. on $K \cap \partial \mathcal{O}$ and $\nabla a^\mathcal{O}=0$ $\mathrm{Cap}_{1,2}-a.e.$ on $K\cap \partial \mathcal{O}$. This means that one can find at least a point $y_0\in K$ such that both $a^\mathcal{O}(y_0)=0$ and $\nabla a^\mathcal{O}(y_0)=0$. Since $a^\mathcal{O}$ is affine, this implies $a^\mathcal{O}=0$, and since this holds for all such connected component $\mathcal{O}$ we deduce that $\phi^{(k+1)}=\phi^{(k)}$ on $\B_k$. 
    
Letting $\phi(x)=\phi^{(k)}(x)$ for $x\in\B_k$, this achieves the construction of $\phi$ globally  on $\R^2$.
\end{proof}

\section{A new monotonicity formula for biharmonic functions}\label{sec_mono}

We consider compact subsets verifying the following hypothesis.
\begin{definition}\label{def_admissible_compact}
A relatively closed set $K\subset \B_1$ is said to be admissible if the set $K\cap\partial\B_r$ is non-empty for every $r\in (0,1)$.
\end{definition}
We will later apply this to sets $K$ such that $0\in K$ and $K\cup\partial\B_1$ is connected, which implies that $K$ is admissible.

We remind the following classical monotonicity formula for harmonic functions: assume $K$ is admissible and $\phi\in H^1_{0,K}(\B_1)$ is harmonic in $\B_1\setminus K$, then
\begin{equation}
\begin{split}\label{eq_lowerorder_monotonicity}
\frac{d}{dr}\left(\frac{1}{r}\int_{ \B_r}|\nabla\phi|^2\right)&=\frac{1}{r}\int_{\partial \B_r}|\nabla\phi|^2-\frac{1}{r^2}\int_{ \B_r}|\nabla \phi|^2\\
&=\frac{1}{r^3}\int_{\partial \B_r}\left(r\partial_r\phi-\frac{1}{2}\phi\right)^2+\left((\partial_\theta \phi)^2-\frac{1}{4}\phi^2\right)
\end{split}
\end{equation}
is nonnegative, and it is zero for every $r$ if and only if there exists $(a,\theta_0)\in\R^2$, such that
$$\phi(re^{i\theta})=ar^\frac{1}{2}\left|\sin\left(\frac{\theta-\theta_0}{2}\right)\right|$$
In this section, we give a higher order version of this formula which, to our knowledge, is new. This formula only relies on the biharmonic equation and the hypothesis on $K$; the minimization of the Griffith energy does not come into play.

\begin{definition}
For any $\phi\in H^2(\B_1)$, for almost every $r\in (0,1)$, we let
\begin{equation}\label{eq_defM}
\mathbf{M}(\phi,r)=\frac{1}{r}\int_{ \B_r}|\nabla^2 \phi|^2+\frac{1}{r}\int_{\partial \B_r}\left(\frac{(\partial_\theta \phi)^2}{r^3}-\frac{(\partial_r \phi)^2}{r}-\frac{\phi\partial_{r,r}\phi}{r}+\frac{5}{2}\frac{\phi\partial_r \phi}{r^2}-\frac{9}{4}\frac{\phi^2}{r^3}\right).
\end{equation}
\end{definition}
Note that the boundary term is a priori well-defined for almost every $r\in (0,1)$. The monotonicity formula below ensures that when $\phi$ is biharmonic on its support, then $\mathbf{M}(\phi,\cdot)$ is absolutely continuous with respect to $r$.
\subsection{Monotonicity for biharmonic functions vanishing on a large set}

The main result of this section is the following.  
\begin{theorem}\label{th_monotonie_disk}
Let $K\subset \B_1$ be an admissible set in the sense of definition \ref{def_admissible_compact}. Let $\phi\in H^{2}_{0,K}(\B_1)$ such that $\Delta^2\phi=0$ in $\B_1\setminus K$. Then the function
$$r\mapsto \mathbf{M}(\phi,r)$$
has the following properties:
\begin{itemize}
\item There exists some universal constant $C>0$ such that 
$$0\leq \mathbf{M}(\phi,r)\leq \frac{C}{r}\int_{\B_r}|\nabla^2\phi|^2$$
for all $r\in (0,1]$.
\item $r\mapsto \mathbf{M}(\phi,r)$ is locally absolutely continuous, non-decreasing, and its weak derivative is
\begin{align*}
\mathbf{M}'(\phi,r)&=\frac{1}{r}\int_{\partial\B_r}\left(\left(\partial_{r,r}\phi-2\frac{\partial_r\phi}{r}+\frac{9}{4}\frac{\phi}{r^2}\right)^2+\frac{2}{r^2}\left(\partial_{r,\theta}\phi-\frac{3}{2}\frac{\partial_\theta\phi}{r}\right)^2+\frac{1}{2r^2}\left(\partial_{r}\phi-\frac{3}{2}\frac{\phi}{r}\right)^2\right)\\
&+\frac{1}{r^5}\int_{\partial\B_r}\left((\partial_{\theta,\theta}\phi)^2-\frac{5}{2}(\partial_{\theta}\phi)^2+\frac{9}{16}\phi^2\right)
\end{align*}
\item $r\mapsto \mathbf{M}(\phi,r)$ is constant if and only if there exists $(a,\theta_0,\theta_1)\in\R^3$ such that
\begin{equation}\label{eq_description32homogeneous}
\forall r\in (0,1),\ \theta\in \frac{\R}{2\pi\Z},\ \phi\left(re^{i\theta}\right)=ar^\frac{3}{2}\sin\left(\frac{\theta-\theta_0}{2}\right)^2\sin\left(\frac{\theta-\theta_1}{2}\right).
\end{equation}
In this case, we have for all $r\in (0,1)$:
$$\mathbf{M}(\phi,r)=\frac{\pi a^2}{4}\left(2+\cos(\theta_0-\theta_1)\right),$$
and either $a=0$ or $K=[0,e^{i\theta_0}]$ up to a $\mathcal{H}^1$-negligible set.
\end{itemize}
\end{theorem}
Note that the function $\sin\left(\frac{\theta-\theta_1}{2}\right)$ is understood to be defined up to a sign, for $\theta\in\frac{\R}{2\pi\Z}\setminus\{\theta_0\}$.
\begin{remark}
This has the following consequences on blow-ups of biharmonic functions: consider $K\subset \B_1$ some connected relatively closed set that joins the origin and $\partial\B_1$, and $\phi\in H^{2}_{0,K}(\B_1)$ such that $\Delta^2\phi=0$ in $\B_1\setminus K$. By \cite[Prop. 6.1]{BCL1}, for any sequence $r_n\to 0$, there is some subsequence $(n_k)_{k\in\N}$ such that $\phi_{r_{n_k}}:=r_{n_k}^{-\frac{3}{2}}\phi(r_{n_k}\cdot)$ converges in $H^2_{\loc}(\R^2)$ to some limit $\psi$ that is biharmonic on $\R^2\setminus L$, where $L$ is a local Hausdorff limit of $K/r_{n_k}$. By the above result, $\psi$ must be of the form
$$\psi(re^{i\theta})=ar^\frac{3}{2}\sin\left(\frac{\theta-\theta_0}{2}\right)^2\sin\left(\frac{\theta-\theta_1}{2}\right)$$
for some $a,\theta_0,\theta_1$ that are constrained by the energy condition
$$\frac{\pi a^2}{4}\left(2+\cos(\theta_0-\theta_1)\right)=\lim_{r\to 0}\mathbf{M}(\phi,r)$$
and $L=e^{i\theta_0}\R_+$ when $a\neq 0$. This does not however guarantee the uniqueness of blow-ups, which is not expected in general. This is further explored in Theorem \ref{th_existence_SIF} in the next section.
\end{remark}
We will apply Theorem \ref{th_monotonie_disk} to the Airy function associated to Griffith minimizers $(u,K)$ in $\B_1$ for which the fracture $K$ connects the origin to $\partial \B_1$. However, since the Airy function is defined implicitly from the displacement field $u$ (by Proposition \ref{prop_def_airy}) and the boundary term in the definition of $\mathbf{M}(\phi,r)$ also involve $\phi$ and its first order derivatives, there seems to be no clean expression of $\mathbf{M}(\phi,r)$ directly in terms of $u$.

In what follows, for the sake of clarity we first prove Theorem \ref{th_monotonie_disk} in conformal polar coordinates. We denote
$$\Cy_\tau=[\tau,+\infty[\times \frac{\R}{2\pi\Z}$$
so that $(t,\theta)\in \Cy_\tau\mapsto e^{-t+i\theta}\in \B_{e^{-\tau}}$ is a conformal map.

We start by giving some technical and computational lemmas, that are consequences of the hypothesis on $K$.
\begin{lemma}\label{lem_bound_origin}
Let $K\subset\B_1$ be an admissible set such that, $\phi\in H^{2}_{0,K}(\B_1)$, then we have
$$\int_{\B_1\setminus \B_{\frac{1}{2}}}\left(|\nabla\phi|^2+\phi^2\right)\lesssim \int_{\B_1\setminus \B_{\frac{1}{2}}}|\nabla^2\phi|^2$$
and
$$\int_{\B_1}\left(\frac{|\nabla\phi|^2}{|x|^2}+\frac{\phi^2}{|x|^4}\right)\lesssim \int_{\B_1}|\nabla^2\phi|^2.$$
\end{lemma}
This guarantees that, letting $\psi(t,\theta):=e^{\frac{3}{2}t}\phi(e^{-t+i\theta})$, the function $\psi$ belongs to $H^2(\Cy_0,e^{-t}dtd\theta)$ with $\Vert\psi\Vert_{H^2(\Cy_0,e^{-t}dtd\theta)}\lesssim \Vert\nabla^2\phi\Vert_{L^2(\B_1)}$. We will first prove the monotonicity in this system of coordinates, where the equality case is clearer.

\begin{proof}
By density it is sufficient to prove this for a function $\phi\in\mathcal{C}^\infty(\overline{\B_1}\setminus K)$ that vanishes in a neighbourhood of $K$.
For any $r\in (0,1)$, by our assumption on $K$ there exists $\theta(r)\in \frac{\R}{2\pi\Z}$ such that $re^{i\theta(r)}\in K$. Then we have for any $\theta\in\frac{\R}{2\pi\Z}$:
\begin{align*}
\phi(re^{i\theta})=\int_{\theta(r)}^{\theta}\partial_{\vartheta}\phi(re^{i\vartheta})d\vartheta,\ \partial_\theta \phi(re^{i\theta})=\int_{\theta(r)}^{\theta}\partial_{\vartheta,\vartheta}\phi(re^{i\vartheta})d\vartheta,\ \partial_r \phi(re^{i\theta})=\int_{\theta(r)}^{\theta}\partial_{r,\vartheta}\phi(re^{i\vartheta})d\vartheta
\end{align*}
The first two estimates follow by integrating this in $r\in\left(\frac{1}{2},1\right)$ or $\left(0,1\right)$.
\end{proof}

\begin{proposition}\label{prop_change_variable_G}
Let $\phi\in H^{2}(\B_1)$, and for any $(t,\theta)\in \Cy_0$, let
$$\psi(t,\theta):=e^{\frac{3}{2}t}\phi(e^{-t+i\theta}).$$
Assume that $\psi\in H^2(\Cy_0,e^{-t}dtd\theta)$. Then for every $\tau\geq 0$, we have
$$\frac{1}{e^{-\tau}}\int_{ \B_{e^{-\tau}}}|\nabla^2 \phi|^2=\mathcal{H}(\psi,\tau)+\int_{\partial\Cy_{\tau}}\Big\{2(\partial_t \psi)^2+\frac{9}{4}\psi^2-\frac{3}{2}\psi\partial_t \psi\Big\}d\theta,$$
where
$$\mathcal{H}(\psi,\tau):=\int_{\Cy_\tau}e^{-(t-\tau)}\Big\{(\partial_{t,t}\psi)^2+2(\partial_{t,\theta}\psi)^2+\frac{3}{2}(\partial_t\psi)^2+(\partial_{\theta,\theta}\psi)^2-\frac{5}{2}(\partial_\theta \psi)^2+\frac{9}{16}\psi^2\Big\}dtd\theta.$$
\end{proposition}
\begin{proof}
For any $r\in (0,1)$, $\theta\in\frac{\R}{2\pi\Z}$, we have
$$\phi(re^{i\theta})=r^\frac{3}{2}\psi(-\log(r),\theta).$$
We compute
\begin{align*}
\partial_\theta\phi(re^{i\theta})&=r^\frac{3}{2}\partial_{\theta}\psi(-\log(r),\theta)\\
\partial_r\phi(re^{i\theta})&=r^\frac{1}{2}\left(\frac{3}{2}\psi -\partial_t\psi\right)(-\log(r),\theta)\\
\partial_{r,\theta}\phi(re^{i\theta})&=r^\frac{1}{2}\left(\frac{3}{2}\partial_{\theta}\psi -\partial_{t,\theta}\psi\right)(-\log(r),\theta)\\
\partial_{r,r}\phi(re^{i\theta})&=r^{-\frac{1}{2}}\left(\frac{3}{4}\psi -2\partial_t\psi+\partial_{t,t}\psi\right)(-\log(r),\theta).
\end{align*}
Then, letting $r=:e^{-t}$, the full Hessian is
\begin{align*}
r|\nabla^2\phi(re^{i\theta})|^2&=r(\partial_{r,r}\phi)^2+r\left(\frac{1}{r}\partial_r\phi+\frac{1}{r^2}\partial_{\theta,\theta}\phi\right)^2+2r\left(\frac{1}{r}\partial_{r,\theta}\phi-\frac{1}{r^2}\partial_\theta\phi\right)^2\\
&=(\partial_{t,t}\psi)^2+2(\partial_{t,\theta}\psi)^2+5(\partial_t\psi)^2+(\partial_{\theta,\theta}\psi)^2+\frac{1}{2}(\partial_\theta\psi)^2+\frac{45}{16}\psi^2\\
&-4\partial_t\psi\partial_{t,t}\psi-6\psi\partial_t\psi+\frac{3}{2}\psi\partial_{t,t}\psi-2\partial_{\theta,\theta}\psi\partial_{t}\psi+3\psi\partial_{\theta,\theta}\psi-2\partial_{\theta}\psi\partial_{t,\theta}\psi,
\end{align*}
where the right-hand side above is evaluated at $(t,\theta)=(-\log(r),\theta)$. We integrate this on $ \B_{e^{-\tau}}$:
\begin{align*}
\int_{ \B_{e^{-\tau}}}|\nabla^2\phi|^2&=\int_{\Cy_\tau}e^{-t}\left((\partial_{t,t}\psi)^2+2(\partial_{t,\theta}\psi)^2+5(\partial_t\psi)^2+(\partial_{\theta,\theta}\psi)^2+\frac{1}{2}(\partial_\theta\psi)^2+\frac{45}{16}\psi^2\right)dt d\theta\\
&+\int_{\Cy_\tau}e^{-t}\left(-4\partial_t\psi\partial_{t,t}\psi-6\psi\partial_t\psi+\frac{3}{2}\psi\partial_{t,t}\psi-2\partial_{\theta,\theta}\psi\partial_{t}\psi+3\psi\partial_{\theta,\theta}\psi-2\partial_{\theta}\psi\partial_{t,\theta}\psi\right)dt d\theta\\
&=\int_{\Cy_\tau}e^{-t}\left((\partial_{t,t}\psi)^2+2(\partial_{t,\theta}\psi)^2+5(\partial_t\psi)^2+(\partial_{\theta,\theta}\psi)^2-\frac{5}{2}(\partial_\theta\psi)^2+\frac{45}{16}\psi^2\right)dt d\theta\\
&+\int_{\Cy_\tau}e^{-t}\left(-4\partial_t\psi\partial_{t,t}\psi-6\psi\partial_t\psi+\frac{3}{2}\psi\partial_{t,t}\psi\right)dt d\theta.
\end{align*}
By integration by part, we simplify the last line above as
\begin{align*}
\int_{\Cy_\tau}e^{-t}\left(-4\partial_t\psi\partial_{t,t}\psi-6\psi\partial_t\psi+\frac{3}{2}\psi\partial_{t,t}\psi\right)dt d\theta&=\int_{\Cy_\tau}e^{-t}\left(-\frac{7}{2}(\partial_t\psi)^2-\frac{9}{4}\psi^2\right)dt d\theta\\
&+e^{-\tau}\int_{\partial\Cy_\tau}\left(2(\partial_t\psi)^2+\frac{9}{4}\psi^2-\frac{3}{2}\psi\partial_t\psi\right) d\theta.
\end{align*}
We obtain
\begin{align*}
\int_{ \B_{e^{-\tau}}}|\nabla^2\phi|^2&=\int_{\Cy_\tau}e^{-t}\left((\partial_{t,t}\psi)^2+2(\partial_{t,\theta}\psi)^2+\frac{3}{2}(\partial_t\psi)^2+(\partial_{\theta,\theta}\psi)^2-\frac{5}{2}(\partial_\theta\psi)^2+\frac{9}{16}\psi^2\right)dt d\theta\\
&+e^{-\tau}\int_{\partial\Cy_\tau}\left(2(\partial_t\psi)^2+\frac{9}{4}\psi^2-\frac{3}{2}\psi\partial_t\psi\right) d\theta,
\end{align*}
which is the result.
\end{proof}

\begin{lemma}\label{lem_boundaryenergy}
Let $K\subset\B_1$ be an admissible set. Let $\phi\in H^{2}_{0,K}(\B_1)$ such that $\Delta^2 \phi=0$ on $\B_1\setminus K$. Let $$\psi(t,\theta):=e^{\frac{3}{2}t}\phi(e^{-t+i\theta}).$$
Assume that $\psi\in H^2(\Cy_0,e^{-t}dtd\theta)$. Let $\zeta\in\mathcal{C}^\infty(\R_{>0},\R)$ such that $\zeta=0$ in a neighbourhood of $0$ and $\zeta'$ has compact support. Then
\begin{align*}
&\int_{\Cy_0}e^{-t}\zeta(t)\Big\{(\partial_{t,t}\psi)^2+2(\partial_{t,\theta}\psi)^2+\frac{3}{2}(\partial_t\psi)^2+(\partial_{\theta,\theta}\psi)^2-\frac{5}{2}(\partial_\theta \psi)^2+\frac{9}{16}\psi^2\Big\}dtd\theta\\
&=-\int_{\Cy_0}e^{-t}\left\{\zeta''(t)\psi\partial_{t,t}\psi+\zeta'(t)\left(2\partial_{t}\psi\partial_{t,t}\psi+2\partial_{\theta}\psi\partial_{t,\theta}\psi+\frac{3}{2}\psi\partial_t \psi\right)\right\}dtd\theta\\
\end{align*}
\end{lemma}

\begin{remark}
If $\partial K$ is sufficiently smooth, so that $\phi$ is $\mathcal{C}^3$ up to the boundary, then the identity above reduces to 
$$\mathcal{H}(\psi,\tau)=\int_{\partial\Cy_\tau}\Big\{\psi\partial_{t,t,t}\psi-\psi\partial_{t,t}\psi-\partial_{t,t}\psi\partial_{t}\psi-2\partial_{t,\theta}\psi\partial_\theta \psi-\frac{3}{2}\psi\partial_t \psi  \Big\}d\theta.$$
as $\zeta$ approaches an indicator function $\chi_{[\tau,+\infty)}$. Note the $3$rd order term, which is not well-defined for a non-smooth set.
\end{remark}

\begin{proof}
$\phi$ is biharmonic in its support, so it is a minimizer of the energy
$$\varphi\in\phi+H^2_0(\B_1\setminus K)\mapsto \int_{\B_1}|\nabla^2 \varphi|^2$$
Let $\eta\in\mathcal{C}^\infty_c(\B_1)$ defined by $\eta(re^{i\theta})=\zeta(-\log(r))$, then by the energy minimality of $\phi$ we have
$$\int_{\B_1}|\nabla^2 \phi|^2\leq\int_{\B_1}|\nabla^2 (e^{\eps\eta}\phi)|^2$$
for every $\eps\in\R$. By Proposition \ref{prop_change_variable_G}, this implies
$$\mathcal{H}(\psi,0)\leq \mathcal{H}(e^{\eps\zeta}\psi,0)$$
for every $\eps\in\R$. Taking the derivative at $\eps=0$, we obtain the relation
$$\int_{\Cy_0}e^{-t}\Big\{\partial_{t,t}\psi\partial_{t,t}(\zeta\psi)+2\partial_{t,\theta}\psi\partial_{t,\theta}(\zeta\psi)+\frac{3}{2}\partial_t\psi\partial_t(\zeta\psi)+\zeta(\partial_{\theta,\theta}\psi)^2-\frac{5}{2}\zeta(\partial_\theta \psi)^2+\frac{9}{16}\zeta\psi^2\Big\}dtd\theta=0,$$
which is (after developing the derivatives) the desired relation.
\end{proof}
\begin{definition}\label{def_W}
For any $\psi\in H^2(\Cy_0,e^{-t}dtd\theta)$, for almost every $\tau\geq 0$, we let
$$\mathcal{M}(\psi,\tau):=\mathcal{H}(\psi,\tau)+\int_{\partial\Cy_\tau}\Big\{(\partial_t \psi)^2-\psi\partial_{t,t}\psi+\psi\partial_t \psi+(\partial_\theta \psi)^2+\frac{3}{4}\psi^2\Big\}d\theta.$$
\end{definition}

\begin{proposition}\label{prop_derivativeW}
Let $K\subset\B_1$ be an admissible set, and $\phi\in H^{2}_{0,K}(\B_1)$ such that $\Delta^2 \phi=0$ in $\B_1\setminus K$. Let $$\psi(t,\theta):=e^{\frac{3}{2}t}\phi(e^{-t+i\theta}).$$
Assume that $\psi\in H^2(\Cy_0,e^{-t}dtd\theta)$. Then the function $\tau\mapsto \mathcal{M}(\psi,\tau)$ is locally absolutely continuous, and its weak derivative is
\begin{equation}\label{eq_derivative_polar_W}
\mathcal{M}'(\psi,\tau)=-\int_{\partial\Cy_\tau}\Big\{(\partial_{t,t}\psi)^2+2(\partial_{t,\theta}\psi)^2+\frac{1}{2}(\partial_{t}\psi)^2+(\partial_{\theta,\theta}\psi)^2-\frac{5}{2}(\partial_\theta \psi)^2+\frac{9}{16}\psi^2\Big\}d\theta.
\end{equation}
\end{proposition}
\begin{proof}
In what follows, we let
$$h(\psi,t):=\int_{\partial\Cy_t}\Big\{(\partial_{t,t}\psi)^2+2(\partial_{t,\theta}\psi)^2+\frac{3}{2}(\partial_t\psi)^2+(\partial_{\theta,\theta}\psi)^2-\frac{5}{2}(\partial_\theta \psi)^2+\frac{9}{16}\psi^2\Big\},$$
so that $\mathcal{H}(\psi,\tau)=\int_{\tau}^{+\infty}e^{\tau-t}h(\psi,t)dt$.

Let $\zeta\in\mathcal{C}^\infty_c(\R_{>0},\R)$, then \eqref{eq_derivative_polar_W} is equivalent to proving that
$$\int_{0}^{+\infty}\zeta'(\tau)\mathcal{M}(\psi,\tau)d\tau=\int_{0}^{+\infty}\zeta(t)h(\psi,t)dt-\int_{\Cy_0}\zeta(t)(\partial_t \psi)^2dtd\theta.$$
Indeed, this is exactly the characterization of the weak derivative of $\mathcal{M}(\psi,\tau)$. We compute directly:
\begin{equation}\label{eq1_proof_prop_derivativeW}
\int_{0}^{+\infty}\zeta'(\tau)\mathcal{M}(\psi,\tau)d\tau=\int_{0}^{+\infty}\zeta'(\tau)\mathcal{H}(\psi,\tau)d\tau+\int_{\Cy_0}\zeta'(t)\Big\{(\partial_t \psi)^2-\psi\partial_{t,t}\psi+\psi\partial_t \psi+(\partial_\theta \psi)^2+\frac{3}{4}\psi^2\Big\}dtd\theta.
\end{equation}
The first term of \eqref{eq1_proof_prop_derivativeW} is simplified as follows:
\begin{align*}
\int_{0}^{+\infty}\zeta'(\tau)\mathcal{H}(\psi,\tau)d\tau&=-\int_{0}^{+\infty}\zeta(\tau)\mathcal{H}'(\psi,\tau)d\tau=\int_{0}^{+\infty}\zeta(\tau)\left(h(\psi,\tau) -\mathcal{H}(\psi,\tau)\right)d\tau\\
&=\int_{0}^{+\infty}\zeta(t)h(\psi,t)dt-\int_{0}^{+\infty}e^{-t}\eta(t)h(\psi,t)dt,
\end{align*}
where $\eta(t):=\int_{0}^{t}e^{s}\zeta(s)ds$. By lemma \ref{lem_boundaryenergy} applied to $\eta$ (note that $\eta$ verifies the correct hypothesis, since $\eta'(t)=e^{t}\zeta(t)$ has compact support), then
\begin{align*}
\int_{0}^{+\infty}\zeta'(\tau)\mathcal{H}(\psi,\tau)d\tau&=\int_{0}^{+\infty}\zeta(t)h(\psi,t)dt+\int_{\Cy_0}e^{-t}\eta''(t)\psi\partial_{t,t}\psi dt d\theta\\
&+\int_{\Cy_0}e^{-t}\eta'(t)\left(2\partial_{t}\psi\partial_{t,t}\psi+2\partial_{\theta}\psi\partial_{t,\theta}\psi+\frac{3}{2}\psi\partial_t \psi\right)dtd\theta\\
&=\int_{0}^{+\infty}\zeta(t)h(\psi,t)dt+\int_{\Cy_0}\zeta'(t)\psi\partial_{t,t}\psi dt d\theta\\
&+\int_{\Cy_0}\zeta(t)\left(\psi\partial_{t,t}\psi+2\partial_{t}\psi\partial_{t,t}\psi+2\partial_{\theta}\psi\partial_{t,\theta}\psi+\frac{3}{2}\psi\partial_t \psi\right)dtd\theta,
\end{align*}
where we replace $e^{-t}\eta'(t)=\zeta(t)$, $e^{-t}\eta''(t)=\zeta(t)+\zeta'(t)$ in the last equality. Now, we simplify the last term by integration by part
\begin{align*}
\int_{0}^{+\infty}\zeta'(\tau)\mathcal{H}(\psi,\tau)d\tau&=\int_{0}^{+\infty}\zeta(t)h(\psi,t)dt+\int_{\Cy_0}\zeta'(t)\psi\partial_{t,t}\psi dt d\theta\\
&-\int_{\Cy_0}\Big\{\zeta(t)(\partial_t \psi)^2+\zeta'(t)\left(\psi\partial_t \psi +(\partial_{t}\psi)^2+(\partial_{\theta}\psi)^2+\frac{3}{4}\psi^2\right)\Big\}dtd\theta.
\end{align*}
We inject this into equation \eqref{eq1_proof_prop_derivativeW}, and with the cancellations we obtain:
\begin{align*}
\int_{0}^{+\infty}\zeta'(\tau)\mathcal{M}(\psi,\tau)d\tau=\int_{0}^{+\infty}\zeta(t)h(\psi,t)dt-\int_{\Cy_0}\zeta(t)(\partial_t\psi)^2,
\end{align*}
which is the result.
\end{proof}
We now study the quadratic form that appears in the derivative $\mathcal{M}'(\psi,\tau)$.
\begin{lemma}\label{lem_valeurpropre}
For any $\psi\in H^2_0([0,2\pi])$, we have
$$\int_{0}^{2\pi}\left((\partial_{\theta,\theta}\psi)^2-\frac{5}{2}(\partial_{\theta}\psi)^2+\frac{9}{16}\psi^2\right)\geq 0,$$
with equality if and only if $\psi\in\mathrm{Span}\left(\sin\left(\frac{\theta}{2}\right)^3,\sin\left(\frac{\theta}{2}\right)^2\cos\left(\frac{\theta}{2}\right)\right)$
\end{lemma}
This is stated in \cite[Lem. 1]{KKLO84} (for the case $\om=2\pi$), however the equality case is not given.
\begin{proof}
We regard $\psi$ as a $2\pi$-periodic function, which still belongs to $H^2(\R/2\pi\Z)$ by the clamping condition. We split $H^2(\R/2\pi\Z)$ into the spaces of even functions and odd functions with respect to $\theta\mapsto -\theta$. Note that this decomposition is orthogonal with respect to the $L^2(\R/2\pi\Z)$, $H^1(\R/2\pi\Z)$ and $H^2(\R/2\pi\Z)$ scalar product, so it is sufficient to prove the inequality in each case.

We now handle separately:
\begin{itemize}
\item Odd functions: in this case, consider $\tilde{\psi}$ the $4\pi$-periodic function defined by
$$\tilde{\psi}(\theta)=\begin{cases}\psi(\theta)&\text{ if }\theta\in [0,2\pi]\\ \psi(-\theta)&\text{ if }\theta\in [-2\pi,0]\end{cases}$$
Again by the clamped condition, $\tilde{\psi}$ belongs to $H^2(\R/4\pi\Z)$ and is even with respect to $\theta\mapsto -\theta$. Thus it may be written as an infinite serie (converging in $H^2(\R/4\pi\Z)$)
$$\tilde{\psi}(\theta)=\sum_{k\in\N}c_k\cos\left(\frac{k\theta}{2}\right)$$
for some sequence $(c_k)_{k\in\N}$ such that $\sum_{k\in\N}k^4|c_k|^2<\infty$, and $c_k=0$ for every even $k$ since $\psi$ is odd. Then
\begin{align*}
\int_{0}^{2\pi}\left((\partial_{\theta,\theta}\psi)^2-\frac{5}{2}(\partial_{\theta}\psi)^2+\frac{9}{16}\psi^2\right)&=\frac{1}{2}\int_{-2\pi}^{2\pi}\left((\partial_{\theta,\theta}\tilde{\psi})^2-\frac{5}{2}(\partial_{\theta}\tilde{\psi})^2+\frac{9}{16}\tilde{\psi}^2\right)\\
&=\pi\sum_{k\in\N}\left(\frac{1}{16}k^4-\frac{5}{8}k^2+\frac{9}{16}\right)|c_k|^2\\
&=\frac{\pi}{16}\sum_{k\in\N}(k^2-1)(k^2-9)|c_k|^2
\end{align*}
This is nonnegative, and vanishes if and only if $\psi\in\mathrm{Span}\left(\cos\left(\frac{\theta}{2}\right),\cos\left(\frac{3\theta}{2}\right)\right)$. The condition $\psi(0)=\psi'(0)=0$ imposes that $\psi$ is a multiple of $\cos\left(\frac{\theta}{2}\right)-\cos\left(\frac{3\theta}{2}\right)$.
\item Even functions: we follow the same lines. Consider $\tilde{\psi}$ the $4\pi$-periodic function defined by
$$\tilde{\psi}(\theta)=\begin{cases}\psi(\theta)&\text{ if }\theta\in [0,2\pi]\\ -\psi(-\theta)&\text{ if }\theta\in [-2\pi,0]\end{cases}$$
$\tilde{\psi}$ belongs to $H^2(\R/4\pi\Z)$ and is odd with respect to $\theta\mapsto -\theta$, thus we may write
$$\tilde{\psi}(\theta)=\sum_{k\in\N}s_k\sin\left(\frac{k\theta}{2}\right)$$
for some sequence $(s_k)_{k\in\N^*}$, with $s_k=0$ for every even $k$ since $\psi$ is even. Then
\begin{align*}
\int_{0}^{2\pi}\left((\partial_{\theta,\theta}\psi)^2-\frac{5}{2}(\partial_{\theta}\psi)^2+\frac{9}{16}\psi^2\right)&=\frac{1}{2}\int_{-2\pi}^{2\pi}\left((\partial_{\theta,\theta}\tilde{\psi})^2-\frac{5}{2}(\partial_{\theta}\tilde{\psi})^2+\frac{9}{16}\tilde{\psi}^2\right)\\
&=\pi\sum_{k\in\N^*}\left(\frac{1}{16}k^4-\frac{5}{8}k^2+\frac{9}{16}\right)|s_k|^2\\
&=\frac{\pi}{16}\sum_{k\in\N^*}(k^2-1)(k^2-9)|s_k|^2
\end{align*}
This is nonnegative, and vanishes if and only if $\psi\in\mathrm{Span}\left(\sin\left(\frac{\theta}{2}\right),\sin\left(\frac{3\theta}{2}\right)\right)$. The condition $\psi(0)=\psi'(0)=0$ imposes that $\psi$ is a multiple of $3\sin\left(\frac{\theta}{2}\right)-\sin\left(\frac{3\theta}{2}\right)$.
\end{itemize}
\end{proof}

\begin{proposition}\label{prop_monotonie_exp}
Let $K\subset\B_1$ be an admissible set. Let $\phi\in H^{2}_{0,K}(\B_1)$ such that $\Delta^2 \phi=0$ on $\B_1\setminus K$. Let $$\psi(t,\theta):=e^{\frac{3}{2}t}\phi(e^{-t+i\theta}).$$
Then $\tau\mapsto \mathcal{M}(\psi,\tau)$ is nonincreasing, and it is constant if and only if there exists $(\theta_0,\theta_1,a)\in\R^3$ such that
$$\phi\left(re^{i\theta}\right)=ar^\frac{3}{2}\sin\left(\frac{\theta-\theta_0}{2}\right)^2\sin\left(\frac{\theta-\theta_1}{2}\right)$$
\end{proposition}
In particular, when $a\neq 0$ this imposes that $K$ must be a segment starting from the origin, up to a set of zero length.
\begin{proof}
Since $K$ is admissible, for every $\tau>0$, there exists $\theta(\tau)\in [0,2\pi)$ such that $e^{-\tau+i\theta(\tau)}\in K$ (otherwise $\partial \B_r$ would separate the origin from $\partial \B_1$), so we have $\psi(\tau,\theta(\tau)+\cdot)\in H^2_0([0,2\pi])$ for almost every $\tau\geq 0$. As a consequence of lemma \ref{lem_valeurpropre}, we have
$$-\mathcal{M}'(\psi,\tau)=\int_{\partial\Cy_\tau}\Big\{\underbrace{(\partial_{t,t}\psi)^2+2(\partial_{t,\theta}\psi)^2+\frac{1}{2}(\partial_{t}\psi)^2}_{\geq 0}\Big\}d\theta+\underbrace{\int_{\partial\Cy_\tau}\Big\{(\partial_{\theta,\theta}\psi)^2-\frac{5}{2}(\partial_\theta \psi)^2+\frac{9}{16}\psi^2\Big\}d\theta}_{\geq 0}.$$
Thus $\mathcal{M}(\psi,\tau)$ is nonincreasing. If it is constant, then by lemma \ref{lem_valeurpropre}, for every $\tau>0$ there exists $a(\tau),\theta_0(\tau),\theta_1(\tau)$ such that
$$\psi(\tau,\theta)=a(\tau)\sin\left(\frac{\theta-\theta_0(\tau)}{2}\right)^2\sin\left(\frac{\theta-\theta_1(\tau)}{2}\right).$$
Moreover, in the case where $\mathcal{M}(\psi,\cdot)$ is constant, we have $\partial_t \psi=0$ almost everywhere: up to a sign choice we identify $\tau\mapsto a(\tau)$ as a constant $a\in\R$. If $a$ is nonzero, then $\theta_0(\tau)$ and $\theta_1(\tau)$ are similarly uniquely defined.
\end{proof}

We may now prove the main result of this section.

\begin{proof}[Proof of Theorem \ref{th_monotonie_disk}]
We split the proof as follows:
\begin{itemize}[label=\textbullet]
\item Definition and monotonicity of $\mathbf{M}$: we let $\psi(t,\theta):=e^{\frac{3}{2}t}\phi(e^{-t+i\theta})$. We first check that for any $\tau\geq 0$,
$$\mathcal{M}(\psi,\tau)=\mathbf{M}(\phi,e^{-\tau}).$$
Indeed, we have
\begin{align*}
\mathcal{M}(\psi,\tau)&=\mathcal{H}(\psi,\tau)+\int_{\partial\Cy_\tau}\Big\{(\partial_t \psi)^2-\psi\partial_{t,t}\psi+\psi\partial_t \psi+(\partial_\theta \psi)^2+\frac{3}{4}\psi^2\Big\}d\theta\text{ by definition of }\mathbf{M}\\
&=\frac{1}{e^{-\tau}}\int_{ \B_{e^{-\tau}}}|\nabla^2 \phi|^2-\int_{\partial\Cy_{\tau}}\Big\{2(\partial_t \psi)^2+\frac{9}{4}\psi^2-\frac{3}{2}\psi\partial_t \psi\Big\}d\theta\\
&+\int_{\partial\Cy_\tau}\Big\{(\partial_t \psi)^2-\psi\partial_{t,t}\psi+\psi\partial_t \psi+(\partial_\theta \psi)^2+\frac{3}{4}\psi^2\Big\}d\theta\text{ by prop. }\ref{prop_change_variable_G}\\
&=\frac{1}{e^{-\tau}}\int_{ \B_{e^{-\tau}}}|\nabla^2 \phi|^2+\int_{\partial\Cy_{\tau}}\Big\{-\psi\partial_{t,t}\psi+\frac{5}{2}\psi\partial_{t}\psi-(\partial_t\psi)^2+(\partial_\theta\psi)^2-\frac{3}{2}\psi^2\Big\}d\theta
\end{align*}
Since $\psi(t,\theta)=e^{\frac{3}{2}t}\phi\left(e^{-t+i\theta}\right)$, we compute, denoting $r:=e^{-t}$,
\begin{align*}
\partial_\theta\psi(t,\theta)&=r^{-\frac{3}{2}}\partial_{\theta}\phi\left(r e^{i\theta}\right)\\
\partial_t\psi(t,\theta)&=\frac{3}{2}r^{-\frac{3}{2}}\phi\left(r e^{i\theta}\right)-r^{-\frac{1}{2}}\partial_{r}\phi\left(r e^{i\theta}\right)\\
\partial_{t,t}\psi(t,\theta)&=\frac{9}{4}r^{-\frac{3}{2}}\phi\left(r e^{i\theta}\right)-2r^{-\frac{1}{2}}\partial_{r}\phi\left(r e^{i\theta}\right)+r^{\frac{1}{2}}\partial_{r,r}\phi\left(r e^{i\theta}\right)
\end{align*}
We simplify
\begin{align*}
&-\psi\partial_{t,t}\psi+\frac{5}{2}\psi\partial_{t}\psi-(\partial_t\psi)^2+(\partial_\theta\psi)^2-\frac{3}{2}\psi^2=\frac{(\partial_\theta \phi)^2}{r^3}-\frac{(\partial_r \phi)^2}{r}-\frac{\phi\partial_{r,r}\phi}{r}+\frac{5}{2}\frac{\phi\partial_r \phi}{r^2}-\frac{9}{4}\frac{\phi^2}{r^3}
\end{align*}
We obtain $\mathcal{M}(\psi,\tau)=\mathbf{M}(\phi,r)$. The expression of $\mathbf{M}'(\phi,r)$ is likewise computed from $\mathcal{M}'(\psi,-\log(r))$.

Finally, by Proposition \ref{prop_monotonie_exp} and lemma $r\in (0,1)\mapsto \mathbf{M}(\phi,r)$ is nondecreasing.
\item Continuity of $\mathbf{M}$: by the trace inequalities 
$$H^2(\B_1,\R)\to H^1(\partial\B_1,\R),$$
the only term of $\mathbf{M}(\phi,r)$ which is not clearly controlled by $\int_{\B_1}|\nabla^2\phi|^2$ is $\int_{\partial\B_1}\phi\partial_{r,r}\phi$. However, using the weak characterization of $\Delta^2\phi=0$ in $\B_1\setminus K$, we have
$$\int_{\B_1}\nabla^2\left(\frac{r^2-1}{2}\phi\right):\nabla^2\phi=\int_{\partial\B_1}\phi\partial_{r,r}\phi,$$
so we conclude that
$$|\mathbf{M}(\phi,1)|\lesssim \int_{\B_1}|\nabla^2\phi|^2$$
By the same proof, we obtain in fact that for any two functions $\phi_1,\phi_2$ satisfying the hypothesis of the theorem, associated to some admissible compact sets $K_1,K_2$, we have
\begin{equation}\label{eq_continuity_M_2}
|\mathbf{M}(\phi_1,1)-\mathbf{M}(\phi_2,1)|\lesssim \Vert \phi_1-\phi_2\Vert_{H^2(\B_1)}\left(\Vert \phi_1\Vert_{H^2(\B_1)}+\Vert \phi_2\Vert_{H^2(\B_1)}\right).
\end{equation}
\item Positivity of $\mathbf{M}$: If $\phi$ vanishes in a neighbourhood of the origin, then $\mathbf{M}(\phi,1)\geq \mathbf{M}(\phi,r)=0$ for some small enough $r$.

Let us now prove the general case by a density argument. For any $\eps>0$, we let $K^{(\eps)}=K\cup\partial\B_\eps$. Let $\eta\in\mathcal{C}^\infty_c(\B_1,\R)$ that equal $1$ in a neighbourhood of the origin, so that $\eta\phi\in H^2_0(\B_1\setminus K,\R)$. Since $K$ contains the origin, by definition of $H^2_{0,K}(\B_1,\R)$, there exists for any $\eps>0$ some function  $\phi^{(\eps,0)}\in \mathcal{C}^\infty_c(\B_1,\setminus K,\R)$ such that $\phi^{(\eps,0)}\equiv 0$ in $\B_\eps$ and
$$\Vert \eta\phi - \phi^{(\eps,0)}\Vert_{H^2(\B_1)}\underset{\eps \to 0}{\longrightarrow}0$$
As a consequence, letting $\phi^{(\eps,1)}=(1-\eta)\phi+\phi^{(\eps,0)}$, we have $\phi^{(\eps,1)}\in H^2_{0,K^{(\eps)}}(\B_1,\R)$ and
\begin{equation}\label{eq_convergence_phieps1}
\Vert \phi - \phi^{(\eps,1)}\Vert_{H^2(\B_1)}\underset{\eps \to 0}{\longrightarrow}0
\end{equation}
Let now $\phi^{(\eps)}\in H^2_{0,K^{(\eps)}}(\B_1,\R)$ be the biharmonic extension of $\phi$ in $\B_1\setminus K^{(\eps)}$, meaning that $\Delta^2 \phi^{(\eps)}=0$ in $\B_1\setminus  K^{(\eps)}$ and $\phi^{(\eps)}-\phi\in H^2_0(\B_1,\R)$. $\phi^{(\eps)}$ vanishes in $\B_\eps$, so the previous point applies: we have
$$\mathbf{M}(\phi^{(\eps)},1)\geq 0.$$
By continuity of $\mathbf{M}(\cdot,1)$ with respect to the $H^2(\B_1)$ norm obtained in equation \ref{eq_continuity_M_2}, it is sufficient to prove that $\phi^{(\eps)}$ converges to $\phi$ in $H^2(\B_1)$. Since $\phi^{(\eps)}-\phi\in H^2_0(\B_1,\R)$, this is verified as follows:
\begin{align*}
\int_{\B_1}|\nabla^2(\phi^{(\eps)}-\phi)|^2&=\int_{\B_1}\left(|\nabla^2\phi^{(\eps)}|^2-|\nabla^2\phi|^2\right)&\text{ by minimality of }\phi\\
&\leq \int_{\B_1}\left(|\nabla^2\phi^{(\eps,1)}|^2-|\nabla^2\phi|^2\right)&\text{ by minimality of }\phi^{(\eps)}\\
&\underset{\eps \to 0}{\longrightarrow}0&\text{ by }\eqref{eq_convergence_phieps1}
\end{align*}
\item Identification of homogeneous solutions: by Proposition \ref{prop_monotonie_exp}, $r\mapsto \mathbf{M}(\phi,r)$ is constant if and only of there exists $a\in\R$, $\theta_0,\theta_1\in\R$ such that
$$\forall r\in (0,1),\ \theta\in\R,\ \phi(re^{i\theta})=ar^\frac{3}{2}\sin\left(\frac{\theta-\theta_0}{2}\right)^2\sin\left(\frac{\theta-\theta_1}{2}\right).$$
In this case, the value of $\mathbf{M}(\phi,r)$ is computed as follows:
\begin{align*}
\mathbf{M}(\phi,1)&=\mathcal{M}(\psi,0)=\mathcal{H}(\psi,0)+\int_{\partial\Cy_0}\Big\{(\partial_t \psi)^2-\psi\partial_{t,t}\psi+\psi\partial_t \psi+(\partial_\theta \psi)^2+\frac{3}{4}\psi^2\Big\}d\theta\\
&=\int_{\partial\Cy_0}\left((\partial_\theta \psi)^2+\frac{3}{4}\psi^2\right)d\theta\text{ because }\mathcal{G}(\psi,0)=0\\
&=\frac{\pi}{4}a^2\left(2+\cos(\theta_0-\theta_1)\right)\text{ after computation.}
\end{align*}
\end{itemize}
\end{proof}

\subsection{Growth of biharmonic functions vanishing on a large set}

The main result of this subsection is a polynomial growth rate for biharmonic functions that vanish on a sufficiently large set, in the sense that no circle arc of $\partial\B_r\setminus K$ may have an angular opening larger than $\om$ for some $\om\in (0,2\pi)$. This result is in fact already contained, with sharp exponent, in the (untranslated) work \cite[Th. 2]{KKLO84}. Our goal is to give an alternative proof based on the monotonicity formula $\mathbf{M}$ which may present independent interest. 

\begin{proposition}\label{prop_poly_growth}
Let $\om\in (0,2\pi)$, then there exists $\kappa_\om,C_\om>0$ such that the following holds. Let $\tau\in (0,1)$, Let $K\subset\B_1$ be an admissible set, such that for every $r\in (\tau,1)$, no connected component of 
$$\left\{\theta\in \frac{\R}{2\pi\Z}:re^{i\theta}\notin K\right\}$$
has length larger than $\om$. Let $\phi\in H^2_{0,K}(\B_1)$, assume that $\Delta^2\phi=0$ in $\B_1\setminus K$. Then
$$\forall r\in (\tau,1),\ \frac{1}{r}\int_{\B_{r}}|\nabla^2\phi|^2\leq C_\om r^{\kappa_\om}\int_{\B_1}|\nabla^2\phi|^2.$$
\end{proposition}

We first remind the following Cacciopoli-type inequality on $\phi$.

\begin{lemma}\label{lem_cacciopoli}
There exists a constant $C>0$ such that for any admissible set $K\subset\B_1$, for any $\phi\in H^2_{0,K}(\B_1)$ such that $\Delta^2\phi=0$ in $\B_1\setminus K$, we have
$$\int_{\B_\frac{1}{2}}|\nabla^2\phi|^2\leq C\int_{\B_1\setminus \B_{\frac{1}{2}}}|\nabla^2 \phi|^2$$
\end{lemma}
\begin{proof}
Consider $\zeta\in\mathcal{C}^\infty_c(\B_1)$ such that $\zeta\equiv 1$ in $ \B_{\frac{1}{2}}$ and $|\nabla^k\zeta|\lesssim 1$ for $k=0,1,2$.  By the energy minimality of $\phi$ we have
$$\forall \eps\in\R,\ \int_{\B_1}|\nabla^2 \phi|^2\leq\int_{\B_1}|\nabla^2 (e^{\eps\zeta^4}\phi)|^2,$$
so
$$\int_{\B_1}\nabla^2(\zeta^4 \phi):\nabla^2\phi=0.$$
Developing $\nabla^2(\zeta^4\phi)$, we obtain
$$\int_{\B_1}\zeta^4|\nabla^2\phi|^2=-4\int_{\B_1}\left(\left(\zeta^3\nabla^2\zeta+3\zeta^2\nabla \zeta\otimes\nabla\zeta\right)\phi:\nabla^2\phi+2\zeta^3(\nabla \zeta\otimes\nabla\phi):\nabla^2\phi\right).$$
By Young inequality, and using the fact that the support of $\nabla\zeta$, $\nabla^2\zeta$ is included in $\B_1\setminus \B_{\frac{1}{2}}$, we obtain that for some universal constant $C>0$:
$$\int_{\B_1}\zeta^4|\nabla^2\phi|^2\leq C\int_{\B_1\setminus \B_{\frac{1}{2}}}(\phi^2+|\nabla\phi|^2)$$
and we conclude with lemma \ref{lem_bound_origin}. 
\end{proof}

\begin{proof}[Proof of Proposition \ref{prop_poly_growth}]
We assume without loss of generality that $\tau\leq \frac{1}{2}$, since otherwise the result is directly true.

By lemma \ref{lem_valeurpropre}, the Rayleigh quotient
$$\inf_{\xi\in H^2_0([0,2\pi])}\frac{\int_{0}^{2\pi}\left((\partial_{\theta,\theta}\xi)^2-\frac{5}{2}(\partial_\theta \xi)^2\right)}{\int_{0}^{2\pi}\xi^2}$$
reaches its minimum with value $-\frac{9}{16}$ for a $2$-dimensional vector space of functions with full support. In particular for any $\om\in (0,2\pi)$, we may find $\delta_\om>0$ such that
$$\inf_{\xi\in H^2_0([0,\om])}\frac{\int_{0}^{2\pi}\left((1-\delta_\om)(\partial_{\theta,\theta}\xi)^2-\frac{5}{2}(\partial_\theta \xi)^2\right)}{\int_{0}^{2\pi}\xi^2}\geq -\frac{9}{16}.$$
We claim that as a consequence, we have
\begin{equation}\label{eq_bound_WH2}
\mathbf{M}'(\phi,r)\gtrsim \frac{\delta_\om}{r}\int_{\partial \B_r}|\nabla^2 \phi|^2
\end{equation}
for any $r\in\left(\tau,1\right)$. Indeed, write $\psi(t,\theta)=e^{\frac{3}{2}t}\phi(e^{-t+i\theta})$, then by Proposition \ref{prop_derivativeW} we have (denoting $t:=-\log(r)$):
\begin{align*}
\mathbf{M}'(\phi,r)&=-e^{t}\mathcal{M}'(\psi,t)\\
&=e^{t}\int_{\partial\Cy_t}\left((\partial_{t,t}\psi)^2+2(\partial_{t,\theta}\psi)^2+\frac{1}{2}(\partial_{t}\psi)^2+(\partial_{\theta,\theta}\psi)^2-\frac{5}{2}(\partial_\theta\psi)^2+\frac{9}{16}\psi^2\right)\\
&\geq e^{t}\int_{\partial\Cy_t}\left((\partial_{t,t}\psi)^2+2(\partial_{t,\theta}\psi)^2+\frac{1}{2}(\partial_{t}\psi)^2+\delta_\om (\partial_{\theta,\theta}\psi)^2\right)\\
&\gtrsim\frac{\delta_\om}{r}\int_{\partial \B_r}|\nabla^2\phi|^2.
\end{align*}
By Theorem \ref{th_monotonie_disk} and lemma \ref{lem_cacciopoli} we have for every $r\in (0,1)$:
$$\mathbf{M}\left(\phi,\frac{r}{2}\right)\lesssim \frac{1}{r}\int_{\B_{\frac{r}{2}}}|\nabla^2\phi|^2\lesssim \frac{1}{r}\int_{\B_r\setminus \B_{\frac{r}{2}}}|\nabla^2\phi|^2.$$
For any $r\in (2\tau,1]$, by integrating \eqref{eq_bound_WH2} from $\frac{r}{2}$ to $r$, we obtain that for some constant $c>0$,
$$c\delta_\om \mathbf{M}\left(\phi,\frac{r}{2}\right)\leq \mathbf{M}\left(\phi,r\right)-\mathbf{M}\left(\phi,\frac{r}{2}\right).$$
This simplifies into $\mathbf{M}\left(\phi,\frac{r}{2}\right)\leq \frac{1}{1+c\delta_\om}\mathbf{M}\left(\phi,r\right)$. By iteration, we obtain
\begin{equation}\label{eq_aux_decrease_M}
\mathbf{M}(\phi,r)\lesssim r^{\kappa_\om}\mathbf{M}(\phi,1).
\end{equation}
for any $r\in [\tau,1]$, where $\kappa_\om:=\frac{\log(1+c\delta_\om)}{\log(2)}$. We conclude, for any $r\in\left[\tau,\frac{1}{2}\right]$
\begin{align*}
\frac{1}{r}\int_{ \B_{r}}|\nabla^2\phi|^2&\lesssim \frac{1}{r}\int_{ \B_{2r}\setminus \B_r}|\nabla^2\phi|^2&\text{ by lemma \ref{lem_cacciopoli}}\\
&\lesssim \delta_\om^{-1}\left(\mathbf{M}(\phi,2r)-\mathbf{M}(\phi,r)\right)&\text{ by eq. }\eqref{eq_bound_WH2}\\
&\leq \delta_\om^{-1}\mathbf{M}(\phi,2r)&\text{ by th. }\ref{th_monotonie_disk}\\
&\leq \delta_\om^{-1}r^{\kappa_\om} \mathbf{M}\left(\phi,1\right)&\text{ by eq. }\eqref{eq_aux_decrease_M}\\
&\lesssim \delta_\om^{-1}r^{\kappa_\om}\int_{\B_1}|\nabla^2\phi|^2&\text{ by lemma }\ref{th_monotonie_disk}
\end{align*}
This prove the result.
\end{proof}

\section{Study of cracktip solutions}\label{sec_stab}

In this section, we study the consequences of the previous monotonicity formula on the solutions of the elasticity equation. In subsection \ref{subsec_stab_cracktip}, we prove a classification of all possible $\frac{1}{2}$-homogeneous minimizers of the Griffith problem : we prove that, up to a sign and a rotation, there is at most one non-zero $\frac{1}{2}$-homogeneous Griffith minimizer, denoted $u_{op}$ (as defined in \eqref{eq_uop}), which corresponds to the ``opening'' of a fracture ; a subtle point is that we do not prove that $u_{op}$ is itself a minimizer.
In subsection \ref{subsec_study_SIF}, we prove the stress intensity factors for a solution of the elasticity equation outside a connected set are always well-defined (but not unique) in Theorem \ref{th_existence_SIF}, and prove a monotonicity of the stress intensity factors through small perturbations in Theorem \ref{th_comparison_SIF} and Corollary \ref{cor_direction_kink}. This subsection is not used in the rest of the paper but may present independent interest.

\subsection{Classification of homogeneous Griffith minimizers}\label{subsec_stab_cracktip}

We start by gathering several known results on the energy release of a crack along a path, for which we refer to \cite{AmestoyLeblond92,ChambolleFrancfortMarigo09,ChambolleFrancfortMarigo10}, see also \cite[Th. 4.1]{AlmiLucardesi18} for a modern reference with explicit constants (note a discrepancy of a factor $2$ in the elastic energy).\bigbreak

Consider some function $u\in LD(\B_1\setminus K)$, where $K=[-1,0]\times\{0\}$, that satisfies the elasticity equation in $\B_1\setminus K$ with homogeneous Neumann condition along $K$. Then $u$ admits the asymptotic
$$u(z)=a+k_{op}u_{op}(z)+k_{sl}u_{sl}(z)+\mathcal{O}(|z|),$$ where $a$ is a rigid motion, $u_{op},u_{sl}$ are $\frac{1}{2}$-homogeneous functions that are respectively the \textit{opening mode} and \textit{sliding mode} (sometimes called \textit{mode I} and \textit{mode II}), and the coordinates $(k_{op},k_{sl})$ are the so-called stress intensity factors (SIF). $u_{op},u_{sl}$ are defined as follows:
\begin{equation}\label{eq_uop}
    u_{op}(re^{i\theta}) = 
    \sqrt{\frac{r}{16  \pi \mu (\lambda + \mu)(\lambda + 2 \mu)}}
    \begin{pmatrix}
        -(\lambda + \mu) \cos(3\theta /2) +(\lambda +5\mu)\cos(\theta/2)\\
        -( \lambda + \mu)\sin(3\theta /2) + (3\lambda +7\mu)\sin(\theta/2)
    \end{pmatrix}
\end{equation}
and
\begin{equation}\label{eq_ush}
    u_{sl}(re^{i\theta}) = 
    \sqrt{\frac{r}{16  \pi \mu (\lambda + \mu)(\lambda + 2 \mu) }}
    \begin{pmatrix}
        (\lambda +\mu) \sin(3\theta /2) +(5\lambda +9\mu) \sin(\theta/2)\\
        -(\lambda +\mu)\cos(3\theta /2) +(\lambda-3\mu)  \cos(\theta/2)
    \end{pmatrix}.
\end{equation}
In both cases $\theta$ is taken between $-\pi$ and $\pi$, with  a fracture at $K$. The prefactor is chosen so that $(u_{op},K)$ and $(u_{sl},K)$ are both stationary for the Griffith energy with respect to shape derivatives.

For good measure, let us also give the expression of the stress tensor $u_{op}$, $u_{sl}$:
$$
\sigma_{op}=\A e(u_{op})=\sqrt{\frac{\mu(\lambda+\mu)}{16\pi (\lambda+2\mu)r}}
\begin{pmatrix}
\cos\left(\frac{5\theta}{2}\right)+3\cos\left(\frac{\theta}{2}\right) & 
\sin\left(\frac{5\theta}{2}\right)-\sin\left(\frac{\theta}{2}\right) \\
\sin\left(\frac{5\theta}{2}\right)-\sin\left(\frac{\theta}{2}\right) &
-\cos\left(\frac{5\theta}{2}\right)+5\cos\left(\frac{\theta}{2}\right) & 
\end{pmatrix}
$$
and
$$
\sigma_{sl}=\A e(u_{sl})=\sqrt{\frac{\mu(\lambda+\mu)}{16\pi (\lambda+2\mu)r}}
\begin{pmatrix}
-\sin\left(\frac{5\theta}{2}\right)-7\sin\left(\frac{\theta}{2}\right) & 
\cos\left(\frac{5\theta}{2}\right)+3\cos\left(\frac{\theta}{2}\right) \\
\cos\left(\frac{5\theta}{2}\right)+3\cos\left(\frac{\theta}{2}\right) &
\sin\left(\frac{5\theta}{2}\right)-\sin\left(\frac{\theta}{2}\right) & 
\end{pmatrix}.
$$
The associated Airy functions are
\begin{equation}\label{eq_phiop}
    \phi_{op}(re^{i\theta}) = \frac{4}{3}\sqrt{\frac{\mu(\lambda+\mu)}{\pi(\lambda+2\mu)}}r^\frac{3}{2}\cos\left(\frac{\theta}{2}\right)^3
\end{equation}
and
\begin{equation}\label{eq_phish}
    \phi_{sl}(re^{i\theta}) = -4\sqrt{\frac{\mu(\lambda+\mu)}{\pi(\lambda+2\mu)}}\sin\left(\frac{\theta}{2}\right)r^\frac{3}{2}\cos\left(\frac{\theta}{2}\right)^2
\end{equation}
and indeed it may be checked that $\mathrm{Co}(\nabla^2\phi_{op}),\mathrm{Co}(\nabla^2\phi_{sl})=\sigma_{op},\sigma_{sl}$.\bigbreak

For any closed set $\Gamma\subset\overline{\B_1}$ such that $\Gamma\cap\partial\B_1=K\cap\partial\B_1$, we may denote by $u^\Gamma$ the minimizer of
$$\int_{\B_1\setminus \Gamma}\A e(v):e(v)$$
among functions $v\in LD(\B_1\setminus\Gamma)$ such that $v|_{\partial\B_1}=u|_{\partial\B_1}$. In particular $u^{K}=u$.

Irwin's formula states that
$$\left.\frac{d}{dt}\right|_{t=0}\int_{\B_1\setminus [-1,t]\times\{0\}}\A e(u^{[-1,t]\times\{0\}}):e(u^{[-1,t]\times\{0\}})=-\left(k_{op}^2+k_{sl}^2\right).$$
Note that with our normalization, there is no prefactor depending on the Lamé coefficients in front (we refer to the constants obtained in \cite[Th. 4.1]{AlmiLucardesi18} - note that the mode I and mode II are switched).

More generally, given some angle $\alpha\in (-\pi,\pi)$, then the stress intensity factors of $u^{K\cup [0,t e^{i\alpha}]}$ at the cracktip $t e^{i\alpha}$ are $(k_{op}^*(\alpha)+o_{t\to 0}(1),k_{sl}^*(\alpha)+o_{t\to 0}(1))$, where
\begin{align*}
k_{op}^*(\alpha)&=F_{1,1}(\alpha)k_{op}+F_{1,2}(\alpha)k_{sl}\\
k_{sl}^*(\alpha)&=F_{2,1}(\alpha)k_{op}+F_{2,2}(\alpha)k_{sl}.
\end{align*}
Here $F(\alpha)\in M_2(\R)$ is a matrix that verifies, by \cite[§A.3]{AmestoyLeblond92}:
\begin{equation}\label{eq_asymptoticF}
F(\alpha)=\begin{pmatrix}1+\mathcal{O}(\alpha^2)& -\frac{3}{2}\alpha+\mathcal{O}(\alpha^2)\\ \frac{1}{2}\alpha+\mathcal{O}(\alpha^2) & 1-\left(\frac{4}{\pi^2}+\frac{3}{8}\right)\alpha^2+\mathcal{O}(\alpha^3)\end{pmatrix}.
\end{equation}
More precise asymptotics are known, but we restrict ourselves to the one we need below. Applying Irwin's formula to the cractip $te^{i\alpha}$ gives the energy release rate in the direction $e^{i\alpha}$:
$$\left.\frac{d}{dt}\right|_{t=0^+}\int_{\B_1\setminus (K\cup [0,t e^{i\alpha}])}\A e(u^{K\cup [0,t e^{i\alpha}]}):e(u^{K\cup [0,t e^{i\alpha}]})=-\left(k_{op}^*(\alpha)^2+k_{sl}^*(\alpha)^2\right).$$

The following proposition is partially contained in \cite[§4]{ChambolleFrancfortMarigo09}, and we report it here for completeness.
\begin{proposition}\label{prop_instability_slide}
Let $K=[-1,0]\times\{0\}$, and let $u\in LD(\B_1\setminus K)$ be $\frac{1}{2}$-homogeneous. Assume $(u,K)$ is a topological minimizer of the Griffith energy, then $u=\pm u_{op}$.
\end{proposition}
\begin{proof}
Since $(u,K)$ is a minimizer, then $u$ verifies the elasticity equation with homogeneous Neumann conditions along $K$. Let $\phi$ be the Airy function of $u$, then $\phi$ is a biharmonic, $\frac{3}{2}$-homogeneous function in $H^2_{0,K}(\B_1)$. Let $g(\theta)=\phi(e^{i\theta})$, then $g$ verifies the following differential equation:
$$
\begin{cases}
g^{(4)}+\frac{5}{2}g^{(2)}+\frac{9}{16}g=0\\
g(-\pi)=g'(-\pi)=g(\pi)=g'(\pi)=0
\end{cases}.
$$
The space of solutions is spanned by $\cos\left(\frac{\theta}{2}\right)^3,\cos\left(\frac{\theta}{2}\right)^2\sin\left(\frac{\theta}{2}\right)$, so $\phi\in\mathrm{Span}(\phi_{op},\phi_{sl})$, meaning $u\in\mathrm{Span}(u_{op},u_{sl})$.

We write $u=k_{op}u_{op}+k_{sl}u_{sl}$ for some stress intensity factors $(k_{op},k_{sl})\neq (\pm 1,0)$. Given any piecewise smooth closed subset $\Gamma\subset\overline{\B_1}$ such that $\Gamma\cap\partial\B_1=K\cap\partial\B_1$, we denote as above $u^{\Gamma}$ to be the minimizer of
$$\int_{\B_1\setminus \Gamma}\A e(v):e(v)$$
among functions $v\in LD(\B_1\setminus \Gamma)$ such that $v|_{\partial\B_1}=u|_{\partial\B_1}$.

Using the energy release rate given above, we have the two relations:
$$
\left.\frac{d}{dt}\right|_{t=0}\mathcal{G}\left(u^{[-1,t]\times\{0\}},[-1,t]\times\{0\}\right)=1-k_{op}^2-k_{sl}^2
$$
and for any $\alpha\in (-\pi,\pi)$,
$$
\left.\frac{d}{dt}\right|_{t=0^+}\mathcal{G}\left(u^{K\cup[0,t e^{i\alpha}]},K\cup[0,t e^{i\alpha}]\right)=1-k_{op}^*(\alpha)^2-k_{sl}^*(\alpha)^2
$$
A necessary condition for $(u,K)$ to be a topological minimizer is thus that
$$k_{op}^2+k_{sl}^2=1$$
and for any $\alpha\in (-\pi,\pi)$, 
$$k_{op}^*(\alpha)^2+k_{sl}^*(\alpha)^2\leq 1$$
Using the asymptotics \eqref{eq_asymptoticF}, we have
\begin{align*}
k_{op}^*(\alpha)^2+k_{sl}^*(\alpha)^2&=k_{op}^2+k_{sl}^2-2\alpha k_{op}k_{sl}+\mathcal{O}(\alpha^2).
\end{align*}
Thus these two condition imply that $k_{op}^2+k_{sl}^2=1$ and $k_{op}k_{sl}=0$. If $(k_{op},k_{sl})=(\pm 1,0)$ we are done, and we now treat the remaining case $(k_{op},k_{sl})=(0,\pm 1)$. In this case the next term in the asymptotic \eqref{eq_asymptoticF} then gives
\begin{align*}
k_{op}^*(\alpha)^2+k_{sl}^*(\alpha)^2&=1+\left(\frac{3}{2}-\frac{8}{\pi^2}\right)\alpha^2+\mathcal{O}(\alpha^3).
\end{align*}
This is strictly larger than $1$ for some small enough, nonzero $\alpha$, which concludes the proof.
\end{proof}

\subsection{Evolution of the stress intensity factors}\label{subsec_study_SIF}

We now state some consequences of the monotonicity formula of Theorem \ref{th_monotonie_disk}, on the blow-ups of solutions to the elasticity equation at cracktips, and on the evolution of the stress intensity factors. Explicitely we have, for any $(k_{op},k_{sl})\in\R^2$:
$$\mathbf{M}(k_{op}\phi_{op}+k_{sl}\phi_{sl},1)=\frac{4\mu(\lambda+\mu)}{3(\lambda+2\mu)}\left(k_{op}^2+3k_{sl}^2\right)$$
Thus the quantity $k_{op}^2+3k_{sl}^2$ can be identified through this monotonicity formula. We first introduce the following object.
\begin{definition}
Let $\Om$ be a bounded open subset of $\R^2$ and $K\subset\Om$ some relatively closed connected subset that is not reduced to a point. We say $u\in\mathrm{El}(\Om;K)$ if $(u,K)$ is an admissible pair such that $u\in LD(\Om\setminus K)$ and
$$\forall v\in LD(\Om\setminus K)\text{ such that }\{v\neq u\}\Subset\Om,\text{ we have }\int_{\Om\setminus K}\A e(u):e(u)\leq \int_{\Om\setminus K}\A e(v):e(v).$$
For any $u\in\mathrm{El}(\Om;K)$, for any $z\in K$, we define $\mathbf{SIF}(u,z)$ to be the set of couples $(k_{op},k_{sl})\in\R^2$ such that
\begin{itemize}
\item When $(k_{op},k_{sl})\neq (0,0)$, we have a convergence of the sequence of pairs 
$$(u_{z,r_n},K_{z,r_n})\underset{n\to+\infty}{\longrightarrow}\left(\mathcal{R}\left(k_{op}u_{op}+k_{sl}u_{sl}\right)\circ\mathcal{R}^{-1},\mathcal{R}(\R_-\times\{0\}\right)$$
for some sequence $r_n\to 0$, and some rotation $\mathcal{R}$.
\item When $(k_{op},k_{sl})=(0,0)$, we have a convergence $e(u_{z,r_n})\chi_{\R^2\setminus K_{z,r_n}}\to 0$ in $L^2_\loc(\R^2,S_2(\R))$ for some sequence $r_n\to 0$.
\end{itemize}
\end{definition}
We could also include the orientation $\mathcal{R}$ of the blow-up in the definition above, but it is not relevant to the discussion below.

\begin{theorem}\label{th_existence_SIF}
Let $\Om$ be a bounded open subset of $\R^2$ and $K\subset\Om$ some relatively closed connected subset that is not reduced to a point. Let $u\in\mathrm{El}(\Om;K)$, and $z\in K$. Then $\mathbf{SIF}(u,z)$ is non-empty, and for any $(k_{op},k_{sl})$ and $(k_{op}',k_{sl}')$ in $\mathrm{SIF}(u,z)$ we have
$$k_{op}^2+3k_{sl}^2=k_{op}'^2+3k_{sl}'^2.$$
\end{theorem}
Note that without further information the quantity $k_{op}^2+k_{sl}^2$ governing the energy release rate may be dependent on the choice of blow-up sequence, due to possible slowly spiralling behavious of $K$.
\begin{proof}
Let $\phi\in H^2_{0,K}$ be the Airy function of $u$. By \cite[Th.2]{KKLO84}, the sequence $(\phi_r)_{r\to 0}$ is uniformly bounded in $H^2(\B_R,\R)$, for every $R>0$. Following \cite[Prop. 6.1]{BCL1}, there exists a subsequence of $(r_n)_{n\in\N}$, that we denote $(r_{n_k})_{k\in\N}$, such that $K_{z,r_{n_k}}$ converges in the locally Hausdorff sense to some closed subset $L$ of $\R^2$, that contains an unbounded component passing through the origin, and $\phi_{z,r_{n_k}}$ converges in $H^2_\loc(\R^2,\R)$ to a limit $\psi\in H^2_{0,L,\loc}(\R^2,\R)$. For every $r>0$, we have
$$\mathbf{M}(\psi,r)=\lim_{k\to +\infty}\mathbf{M}(\phi(z+\cdot),r_{n_k}r)=\lim_{\rho\to 0}\mathbf{M}(\phi(z+\cdot),\rho)$$
which does not depend on $r$, thus by Theorem \ref{th_monotonie_disk}, $\psi$ is $\frac{3}{2}$-homogeneous, and up to a rotation $\mathcal{R}$ of the coordinate (which we assume to be the identity without loss of generality), $\psi=k_{op}\phi_{op}+k_{sl}\psi_{sl}$ for some $(k_{op},k_{sl})\in\R^2$, defined by
$$\frac{4\mu(\lambda+\mu)}{3(\lambda+2\mu)}\left(k_{op}^2+3k_{sl}^2\right)=\mathbf{M}(\psi,1)=\lim_{\rho\to 0}\mathbf{M}(\phi(z+\cdot),\rho)$$
so it does not depend on the sequence $(r_n)_n$. We now differentiate two cases.
\begin{itemize}
\item $(k_{op},k_{sl})\neq (0,0)$. Since $\psi$ is biharmonic in $\R^2\setminus L$, this implies $L=\R_-\times\{0\}$, and by \cite[Th. 6.4]{BCL1}, the pair $(u_{z,r_n},K_{z,r_n})$ converges (in the sense of pairs) to $(k_{op}u_{op}+k_{sl}u_{sl},\R_-\times\{0\})$.
\item $(k_{op},k_{sl})=(0,0)$. Then $\nabla^2\phi_{z,r_{n_k}}$ converges to $0$ in $L^2_\loc(\R^2,S_2(\R))$, and so does $$e(u_{z,r_{n_k}})\chi_{\R^2\setminus K_{z,r_{n_k}}}=\A^{-1}\mathrm{Co}(\nabla^2\phi_{z,r_{n_k}}).$$
\end{itemize}
This ends the proof.
\end{proof}
\begin{theorem}\label{th_comparison_SIF}
Let $K$ a relatively closed connected subset of $\B_1$ that connects the origin to $\partial\B_1$. Let $u\in \mathrm{El}(\B_1,K)$, and let $(\gamma^\eps)_{\eps\to 0}$ be a sequence of compact connected subset of $\B_\eps$ that contain the origin, $u^\eps\in\mathrm{El}(\B_1,K\cup\gamma^\eps)$ such that $(u^\eps-u)|_{\partial\B_1\setminus \overline{K}}=0$, and $z^\eps\in \B_\eps\cap(K\cup\gamma^\eps)\setminus\{0\}$. Finally, let
$$(k_{op},k_{sl})\in\mathbf{SIF}(u,0),\ (k_{op}^\eps,k_{sl}^\eps)\in\mathbf{SIF}(u^\eps,z^\eps).$$
Then
$$\limsup_{\eps \to 0}\left((k_{op}^{\eps})^2+3(k_{sl}^\eps)^2\right)\leq k_{op}^2+3k_{sl}^2$$
and if the equality is reached for some sequence $\eps_n\to 0$, then either $(k_{op},k_{sl})=(0,0)$, or some subsequence of $(K\cup\gamma^{\eps_n})/|z^{\eps_n}|$ converges (in the local Hausdorff sense) to $\mathcal{R}\left[(-\infty,1]\times\{0\}\right]$, for some rotation $\mathcal{R}$.
\end{theorem}

\begin{proof}
We define $\phi,\phi^\eps$ the Airy functions of $u,u^\eps$. $\phi^\eps$ belongs to $H^2_{0,K\cup\gamma^\eps}(\B_1,\R)$ which is a subset of $H^2_{0,K}(\B_1,\R)$, and we consider the following twisted scalar product on $H^2_{0,K}(\B_1,\R)$:
$$b_\A( \psi,\eta)=\int_{\B_1}\A^{-1}\mathrm{Co}(\nabla^2\psi):\mathrm{Co}(\nabla^2\eta).$$
We remind that this is a scalar product by the estimates of lemma \ref{lem_bound_origin}. Let $\psi\in H^2_{0,K\cup\gamma^\eps}(\B_1,\R)$, then
\begin{align*}
\int_{\B_1\setminus K}\A^{-1}\mathrm{Co}(\nabla^2(\phi^\eps-\phi)):\mathrm{Co}(\nabla^2\psi)&=\int_{\B_1\setminus K} e(u^\eps-u):\mathrm{Co}(\nabla^2\psi)=0.
\end{align*}
Thus $\phi^\eps$ is the $b_\A$-orthogonal projection of $\phi\in H^2_{0,K}(\B_1,\R)$ on the subspace $H^2_{0,K\cup\gamma^\eps}(\B_1,\R)$. 

Consider $\chi^\eps\in\mathcal{C}^\infty_c(\B_{2\eps},\R)$, such that $\chi^\eps\equiv 1$ on $\B_{\eps}$ and $|\nabla^k \chi^\eps|\lesssim \eps^{-k}$ for $k=0,1,2$. Then $(1-\chi^\eps)\phi\in H^2_{0,K\cup\gamma^\eps}(\B_1,\R)$, so we can estimate
\begin{align*}
\int_{\B_1}|\nabla^2(\phi-\phi^\eps)|^2&\lesssim \int_{\B_1}|\nabla^2(\chi^\eps\phi)|^2\\
&\lesssim \int_{\B_{2\eps}}|\nabla^2\phi|^2&\text{ by lemma }\ref{lem_bound_origin}\\
&\lesssim \eps\int_{\B_1}|\nabla^2\phi|^2 &\text{ by \cite[Th. 2]{KKLO84}}.
\end{align*}
By definition of $(k_{op},k_{sl})$, and up to a rotation of the space, there exists some sequence $r_n\to 0$ such that $(u_{r_n},K_{r_n})$ converge to $(k_{op}u_{op}+k_{sl}u_{sl},\R_-\times\{0\})$ (or simply $u_{r_n}\to 0$ in the case $(k_{op},k_{sl})=(0,0)$), and $\phi_{r_n}$ converge to $k_{op}\phi_{op}+k_{sl}\phi_{sl}$.

Likewise, for any $\eps$, there is some sequence $r_n^\eps\to 0$ and some rotation $\mathcal{R}^\eps$ such that $(u^\eps_{r_n},\mathcal{R}^\eps(K^\eps_{r_n}))$ converge to $\left(\mathcal{R}^\eps\left(k^\eps_{op}u_{op}+k^\eps_{sl}u_{sl}\right)\circ(\mathcal{R}^\eps)^{-1},\mathcal{R}^\eps(\R_-\times\{0\})\right)$, and $\phi_{r_n}$ converges to $\left(k^\eps_{op}\phi_{op}+k^\eps_{sl}\phi_{sl}\right)\circ(\mathcal{R}^\eps)^{-1}$.

Let $\delta>0$, we consider some sufficiently large $n(\delta)$ such that $\Vert \phi_{r_n}-k_{op}\phi_{op}-k_{sl}\phi_{sl}\Vert_{H^2(\B_1)}\leq \delta$. By the previous estimate, $\Vert \phi_{r_{n(\delta)}}^\eps-\phi_{r_{n(\delta)}}\Vert_{H^2(\B_1)}\lesssim  \sqrt{r_{n(\delta)}^{-1}\eps}\Vert\phi\Vert_{H^2(\B_1)}$, so
$$\Vert \phi^\eps_{r_{n(\delta)}}-k_{op}\phi_{op}-k_{sl}\phi_{sl}\Vert_{H^2(\B_1)}\lesssim \delta+\sqrt{r_{n(\delta)}^{-1}\eps}\Vert\phi\Vert_{H^2(\B_1)}.$$
Thus, for sufficiently $\eps$ such that $|z^\eps|\leq \frac{1}{2}r_{n(\delta)}$, we have
$$\mathbf{M}(\phi^\eps_{z^\eps,\frac{1}{2}r_n},1)\leq \mathbf{M}\left((k_{op}\phi_{op}+k_{sl}\phi_{sl})_{z^\eps,\frac{1}{2}},1\right)+C\left(\delta+\sqrt{r_{n(\delta)}^{-1}\eps}\Vert\phi\Vert_{H^2(\B_1)}\right)$$
for some constant $C\lesssim 1$. Using the monotonicity of $\mathbf{M}$ on the left-hand side, this implies
\begin{align*}
\mathbf{M}(k_{op}^\eps u_{op}+k_{sl}^\eps u_{sl},1)&\leq \mathbf{M}(\phi^\eps_{z^\eps,\frac{1}{2}r_{n(\delta)}},1)\\
&\leq \mathbf{M}\left((k_{op}\phi_{op}+k_{sl}\phi_{sl})_{z^\eps,\frac{1}{2}},1\right)+C\left(\delta+\sqrt{r_{n(\delta)}^{-1}\eps}\Vert\phi\Vert_{H^2(\B_1)}\right)
\end{align*}
Taking $\eps\to 0$, we obtain
\begin{align*}
\limsup_{\eps\to 0}\mathbf{M}(k_{op}^\eps u_{op}+k_{sl}^\eps u_{sl},1)&\leq \mathbf{M}\left(k_{op}\phi_{op}+k_{sl}\phi_{sl},1\right)+C\delta
\end{align*}
This is verified for any $\delta>0$, hence the result.

Assume now that $(k_{op},k_{sl})\neq (0,0)$ and that we are in the equality case for some sequence $\eps_k\to 0$, meaning that
$$\lim_{k\to +\infty}\mathbf{M}(k_{op}^{\eps_k}\phi_{op}+k_{sl}^{\eps_k}\phi_{sl},1)=\mathbf{M}(k_{op}\phi_{op}+k_{sl}\phi_{sl},1).$$
Denote
$$\psi^\eps:=\phi^\eps_{z^{\eps},|z^{\eps}|},$$
we prove that the sequence $(\psi^{\eps_k})_k$ converges, up to a rotation, to a non-zero combination of $\phi_{op},\phi_{sl}$, and this is sufficient to imply the result.

Let $R>0$. For any sufficiently small $\eps$, we have by \cite[Th. 2]{KKLO84} that $\Vert \psi^\eps\Vert_{H^2(\B_R)}\lesssim \Vert \phi^{\eps}\Vert_{H^2(\B_1)}(\lesssim \Vert \phi\Vert_{H^2(\B_1)})$, so following \cite[Prop. 6.1]{BCL1} we may assume that some subsequence of $\psi^{\eps_k}$ converges in $H^2_\loc(\R^2,\R)$ to a limit $\psi$, and $(K\cup\gamma^\eps)_{z^\eps,|z^\eps|}$ converges in the locally Hausdorff sense to some limit $L$, which contains an unbounded component passing through the origin, such that $\psi\in H^2_{0,L,\loc}(\R^2,\R)$ and $\Delta^2\psi=0$ in $\R^2\setminus L$.

For any $r>0$, we have $\mathbf{M}(\psi^\eps,r)\geq\mathbf{M}(k_{op}^\eps+k_{sl}^\eps,1)$, and at the same time
$$\mathbf{M}(\psi^\eps,r)\leq \mathbf{M}(\phi^\eps_{z^\eps,\frac{1}{2}r_{n(\delta)}})\leq \mathbf{M}\left((k_{op}\phi_{op}+k_{sl}\phi_{sl})_{z^\eps,\frac{1}{2}},1\right)+C\left(\delta+\sqrt{r_{n(\delta)}^{-1}\eps}\Vert\phi\Vert_{H^2(\B_1)}\right)$$
so passing to the limit in $\eps_k\to 0$ we get
$$\lim_{k\to +\infty}\mathbf{M}(k_{op}^{\eps_k}\phi_{op}+k_{sl}^{\eps_k}\phi_{sl},1)\leq \mathbf{M}(\psi,r)\leq\mathbf{M}(k_{op}\phi_{op}+k_{sl}\phi_{sl},1)+C\delta$$
This is verified for arbitrarily small $\delta$, thus we obtain $\mathbf{M}(\psi,r)=\mathbf{M}(k_{op}\phi_{op}+k_{sl}\phi_{sl},1)$ for all $r>0$. This implies that $\psi$ is, up to some rotation, a non-zero combination of $\phi_{op},\phi_{sl}$, which ends the proof.
\end{proof}

Taking $K=[-1,0]$, $\gamma^\eps=[0,\eps e^{i\alpha}]$, $z^\eps=\eps e^{i\alpha}$, for some $\alpha\in (-\pi,\pi)$, and $u=k_{op}u_{op}+k_{sl}u_{sl}$, this has the direct corollary.
\begin{corollary}\label{cor_direction_kink}
For any $(k_{op},k_{sl})\in\R^2\setminus\{(0,0)\}$, $\alpha\in (-\pi,\pi)\setminus\{0\}$, denote
\begin{align*}
k_{op}^*(\alpha)&=F_{1,1}(\alpha)k_{op}+F_{1,2}(\alpha)k_{sl}\\
k_{sl}^*(\alpha)&=F_{2,1}(\alpha)k_{op}+F_{2,2}(\alpha)k_{sl},
\end{align*}
Then
$$k_{op}^*(\alpha)^2+3k_{sl}^*(\alpha)^2<k_{op}^2+3k_{sl}^2.$$
\end{corollary}

\begin{proof}
We let $K=[-1,0]\times\{0\}$, $u=k_{op}u_{op}+k_{sl}u_{sl}$ with its associated Airy function $\phi=k_{op}\phi_{op}+k_{sl}\phi_{sl}$. We denote $u^{\eps}:=u^{K\cup [0,\eps e^{i\alpha}]}$, and $\phi^{\eps}$ its associated Airy function. From \cite{AmestoyLeblond92} we have
$$\mathbf{SIF}(u^\eps,\eps e^{i\alpha})=\left\{\left(k_{op}^*(\alpha)+o_{\eps\to 0}(1),k_{sl}^*(\alpha)+o_{\eps\to 0}(1)\right)\right\}$$
So by Theorem \ref{th_comparison_SIF} we directly get
$$k_{op}^*(\alpha)^2+3k_{sl}^*(\alpha)^2\leq k_{op}^2+3k_{sl}^2.$$
Assume that we are in the equality case. By Theorem \ref{th_comparison_SIF}, this implies that $([-1,0]\times \{0\}\cup[0,\eps e^{i\alpha}])_{\eps e^{i\alpha},\eps}$ converges to a half-line for some sequence $\eps_k\to 0$. Since $\alpha$ is non-zero, this limit set is instead a broken line, which is a contradiction.
\end{proof}

This has the following consequences:
\begin{itemize}
\item The best direction for the energy release rate in the pure opening regime $k_{sl}=0$, is $\alpha=0$. Indeed it is discussed in \cite[Rem. 4.5]{ChambolleFrancfortMarigo10} that for a cracktip solution with a pure opening mode, the optimal direction for energy release should be to continue the crack in a straight line (i.e. with no kink). This is supported by numerical evidence and by the asymptotic of $F(\alpha)$ for small $\alpha$, and is known in the scalar case per \cite[Prop. 5.4]{Lemenant2016}, but there is no previous existing proof for a general $\alpha\in (-\pi,\pi)$ to the authors' knowledge. If $k_{sl}=0$, then following the notations above we have
$$\forall \alpha\in (-\pi,\pi)\setminus\{0\},\ k_{op}^*(\alpha)^2+k_{sl}^*(\alpha)^2\leq k_{op}^*(\alpha)^2+3k_{sl}^*(\alpha)^2<k_{op}^2$$
\item For any $(k_{op},k_{sl})\in\R^2\setminus\{(0,0)\}$, and any $\alpha\in (-\pi,\pi)\setminus\{0\}$ such that $k_{op}^*(\alpha)^2+k_{sl}^*(\alpha)^2\geq k_{op}^2+k_{sl}^2$, then $|k_{sl}^*(\alpha)|<|k_{sl}|$ and $|k_{op}^{*}(\alpha)|>|k_{op}|$.
\item For any $(k_{op},k_{sl})\in\R^2\setminus\{(0,0)\}$, and any $\alpha\in (-\pi,\pi)\setminus\{0\}$, we have an explicit bound on the maximal energy release rate before and after an angle: $k_{op}^*(\alpha)^2+k_{sl}^*(\alpha)^2<3\left(k_{op}^2+k_{sl}^2\right)$.
\end{itemize}

\section{Regularity for topological minimizers}\label{sec_epsreg}

 The purpose of this section is to prove that if the crack set of a minimizer is sufficiently flat in a ball, then it is $C^1$-regular in a slightly smaller ball. Such an $\varepsilon$-regularity result was already known under an additional topological assumption, namely a separation condition. Here, we remove this assumption and show that flatness alone is sufficient. The key ingredient is a general geometric property of connected components of minimizers, namely the chord-arc property established in \cite{LLL}. Using this property, we are able to dispense with the separation assumption in the $\varepsilon$-regularity result.

\subsection{The chord-arc property}\label{subsec_chordarc}

The following proposition can be found in \cite{LLL}, which is the Griffith analogue to what was proved before by David in \cite[Proposition 68.16]{d}. For this statement only we adopt the notion of almost minimizers like in \cite{LLL}, although in the present paper we restrict ourselves to the case of minimizers, which corresponds to $h = 0$.

\begin{proposition}[Local-John property] {\rm \cite[Proposition 3.9]{LLL}} \label{localJohn}There are constants $C_J\geq 1$ and $\tau_J>0$ such that whenever $(u,K)$ is a Griffith almost minimizer in $\Omega$ with gauge $h$, $x\in \Omega \setminus K$, and $R\in (0,{\rm dist}(x,\partial \Omega))$ is also such that $h(R)\leq \tau_J$, then there is a $1$-Lipschitz mapping $\gamma: [0,R] \to \Omega \setminus K$ such that $\gamma(0)=x$ and 
$$ {\rm dist}(\gamma(t),K)\geq C_J^{-1}t \quad \text{ for } t\in [0,R].$$
\end{proposition}

From the above proposition we can prove that $K$ is (locally) a chord-arc tree. For the recall, we first give the usual definition of a chord-arc set. 

\begin{definition}[Chord-arc]\label{chordarc}
Let $K$ be a closed subset of $\R^2$ which is path-connected. We say that $K$ is \emph{chord-arc} with constant $C_0>0$ if the following property holds: for every $x,y \in K$, there exists an injective, connected curve $\Gamma \subset K$ containing $x$ and $y$ and satisfying
$\HH^1(\Gamma)\leq C_0 |x-y|$.
\end{definition}

As a matter of fact, it could be proved that a global  Griffith minimizer with $K$ connected, satisfies the chord-arc property. On the other hand,  a Griffith almost-minimizer satisfies a local version of chord-arc that is slightly weaker, and that we call \emph{local chord-arc}.  This weaker notion applies both to global and to local minimizers, and  will be established  below following an approach due to David \cite{d}.

\begin{definition}[Local-chord-arc]\label{localchordarc}
Let $K\subset \Omega$ be a relatively closed subset of $\R^2$. We say that $K$ satisfies the \emph{local-chord-arc} property (or is \emph{locally chord-arc}) with constants $C_0>0$ in $\Omega$ if the following property holds: let $\Gamma\subset K$ be a simple curve in $K$. Assume that $\Gamma\subset \B_{x_0,r}$ for some disk satisfying $\B_{x_0,C_0r}\subset \Omega$. Call $x$ and $y$ the endpoints of $\Gamma$. Then $\mathcal{H}^1(\Gamma)\leq C_0|x-y|$. 
\end{definition}

The proof of the following proposition is inspired by \cite[Exercice 57.29]{d} and the statement is similar to  \cite[Proposition 68. 18]{d}.

\begin{proposition}[Local-chord-arc property] \label{locchard}  Let $C_J,\tau_J$ be the constants of Proposition \ref{localJohn}, and $C_A$ be the constant of upper Ahlfors-regularity of $K$.   Let $(u,K)$ be
a Griffith topological minimizer in $\Omega\subset \R^2$. Then   $K$ is locally chord-arc in $\Omega$ with constant $C_0:=400 C_A C_J$.
\end{proposition}
\begin{remark} For a global minimizer it could be proved using the same strategy of proof as the one of Proposition \ref{locchard}  that   $K$ is locally chord-arc in $\R^2$ with constants $C_0$: in other words, if $K$ is furthermore connected then $K$ is chord-arc as defined in Definition \ref{chordarc}.
\end{remark}

\begin{remark} A similar statement could be obtained for almost-minimizers with the same proof.
\end{remark}

\begin{proof}[Proof of Proposition \ref{locchard}] Let $C_J,\tau_J$ be the constants of Proposition \ref{localJohn},  $C_A$ be the constant of upper Ahlfors-regularity of $K$, and $C_0:=  400 C_A C_J$. Let us recall  that $C_A,C_J\geq 1$.

We argue by contradiction. If Proposition \ref{locchard} is false, then it means that one can find two points $x,y \in K \subset \B_{x_0,r}\subset \B_{x_0,C_0r}\subset \Omega$, and a simple arc $\Gamma \subset K\cap \B_{x_0,r}$ connecting them   with 
\begin{eqnarray}
\mathcal{H}^1(\Gamma)\geq C_0 |x-y|. \label{estimLenght}
\end{eqnarray} In particular it implies that  ${\rm diam}(\Gamma)\geq C_0C_A^{-1}|x-y|/4$ otherwise by Ahlfors-regularity of $K$ we would have $\mathcal{H}^1(\Gamma)\leq C_0|x-y|/2$, which is a contradiction with \eqref{estimLenght}.

We call $\widehat \Gamma $ the loop obtained by adding the segment $[x,y]$ to the curve $\Gamma$. The idea to finish the proof of the Proposition is rather simple, by finding a point in the bounded component of $\R^2\setminus \widehat \Gamma$ relatively far away from $[x,y]$, and  then use an escape path from that point given by Proposition~\ref{localJohn} that must necessarily escape the component by the segment $[x,y]$, reaching a contradiction. But the full argument needs some details and the constants have to be fixed in a proper way.

 So we denote $d:= {\rm diam}(\Gamma)\geq C |x-y|$, with $C=C_0C_A^{-1}/4$ and we let  $z \in \Gamma$ be a point such that $|z-x|\geq d/2$. The next goal is to  find a point $\xi$  in the ball $\B_{z,10^{-2}d}$  that belongs to a bounded connected component of $\R^2 \setminus \widehat \Gamma$. This would be easier if $\widehat \Gamma$ would be Jordan curve, which is not sure in general even if $\Gamma$ was injective, because $\Gamma$ could meet the segment $[x,y]$ at many other points than $x$ and $y$. For that purpose we will modify a little $\widehat \Gamma$ as follows: we start from $z$ and follow the curve $\Gamma$ in the direction of $x$, up to finding a first touching point $a\in [x,y]$. Then we follow $\Gamma$ in the other direction and find another first touching point $b\in [x,y]$. We denote by $\Gamma_1$ the portion of injective curve  from $a$ to $z$ in $\Gamma$, and $\Gamma_2$ the the portion from $z$ to $b$. Then the loop formed by the segment $[a,b]$ together with $\Gamma_1\cup \Gamma_2$, form now a Jordan curve containing $z$, denoted by $\widehat \Gamma'$.  By  Jordan's theorem we know $\R^2\setminus \widehat \Gamma'$ has two connected components, one bounded and one unbounded, and  that $z$ belongs to closure of the bounded one. In other words there exists $\xi \in \B_{z,10^{-2}d}$ that belongs to the bounded connected component of $\R^2 \setminus \widehat \Gamma'$.

Now we use the existence of escape paths given by  Proposition \ref{localJohn}. Let us check that the proposition applies, for a curve starting at $\xi$. We know that $\widehat \Gamma' \subset \B_{x_0,r}$ and $\B_{x_0,8C_Jr}\subset \Omega$. This implies that the bounded connected component enclosed by $\widehat \Gamma'$ must also be strictly contained in $\B_{x_0,r}$ and therefore ${\rm dist}(\xi, \partial \Omega) > 4C_Jr$. This allows to apply Proposition \ref{localJohn} with $R=4C_Jr$, which provides the existence of a   $1$-Lipschitz mapping $\gamma: [0,R] \to \Omega \setminus K$ such that 
\begin{eqnarray}
{\rm dist}(\gamma(t),K)\geq C_J^{-1}t \quad \text{ for } t\in [0,R]. \label{estdi}
\end{eqnarray}
In particular for $t_0=2C_Jr$ we have 
\begin{eqnarray}
{\rm dist}(\gamma(t_0),K)\geq 2r, \notag 
\end{eqnarray}
and since $K \cap \B_{x_0,r}\not = \emptyset$, this means that $\gamma(t_0)\in \Omega \setminus \B_{x_0,r}$, and in particular $\gamma(t_0)$ belongs to the unbounded connected component of $\R^2 \setminus \widehat \Gamma'$. But this means that necessarily, the curve $\gamma$ touches the loop somewhere, in other words $\gamma(t_1)\in \widehat \Gamma'$ for some $t_1\in (0,t_0)$, and since $\gamma(t) \not \in K$ for all $t$, we deduce that $\gamma(t_1)\in [a,b]$.

Let us check that this is not possible. Indeed, by construction we have:
\begin{eqnarray}
 |\xi-\gamma(t_1)|&\geq& |z -x|- |z-\xi| - |x-\gamma(t_1)| \notag \\
 &\geq&  d/2-10^{-2}d-C^{-1}d\geq   d/3,
 \end{eqnarray}
 because $C^{-1}=4C_A/ C_0\leq 10^{-2}$ by definition of $C_0$.  And since $\gamma$ is $1$-Lipschitz and $\xi=\gamma(0)$ we deduce that $t_1\geq |\gamma(0)-\gamma(t_1)|\geq d/3$. But then \eqref{estdi} says that 
 $${\rm dist}(\gamma(t_1),K)\geq C_J^{-1}t_1\geq C_J^{-1}d/3.$$
 On the other hand, $\gamma(t_1)$ belongs to $[x,y]$ with $x,y\in K$ and by assumption $d\geq C |x-y|$. This means that ${\rm dist}(\gamma(t_1),K)\leq C^{-1}d$, or in other words $C_0\leq 12C_JC_A$ which is a contradiction with the definition of $C_0$. This achieves the proof of the Proposition.
\end{proof}

\subsection{Flatness implies separation}

\begin{definition}
Let $K$ be a closed subset of $\R^2$,  $x\in \R^2$ and $r>0$ be such that $\beta_K(x,r) \leq 1/2$ and $r<{\rm diam}(K)$. 
We say that $K$ \emph{separates}  in $\B_{x,r}$  if the balls $D^\pm_{x,r}=\B_{x\pm \frac{3}{4}\nu r,\frac{1}{4}r}$, where $\nu$ is a unit vector such that the $\beta_K(x,r)$ is reached for the affine line $\{z:(z-x)\cdot\nu=0\}$, are contained into two different connected components of $\overline{\B_{x,r}}\setminus K$.
\end{definition}

We will frequently use the following geometrical fact: if $\beta_K(x,r)\leq \frac{\rho}{2}$ for some $\rho\in [0,1]$, and $K$ separates $D^\pm_{x,r}$ in $\B_{x,r}$, then $K$ separates $D^{\pm}_{x,\rho r}$ in $\B_{x,\rho r}$.\bigbreak

The aim of this section is to prove a topological fact about  chord-arc sets that will be needed later in view of establishing an  epsilon-regularity property for connected Griffith minimizers.

\begin{proposition} \label{propositionSepare}
Let $K\subset \Omega  \subset  \R^2$ be a relatively closed, connected  subset of $\Om$ such that $\mathcal{H}^1(K)<+\infty$, assume that $K$ is locally chord-arc in $\Omega$ with constant $C_0$. Let $x_0\in K$, then there exists $R>0$ (depending on $K,\Om,x_0$) such that for any $x\in \B_{x_0,R}\cap K$, for any $r>0$ such that $\B_{x,r}\subset\B_{x_0,R}$, we have the following implication: if 
$$\beta_K(x,r)\leq \frac{1}{8\left(2+(1+\sqrt{5})C_0\right)},$$
then $K$ separates $D^\pm_{x,\frac{1}{2}r}$ in $\B_{x,\frac{1}{2}r}$.
\end{proposition}
\begin{proof}
We write $\eta=\frac{1}{8\left(2+(1+\sqrt{5})C_0\right)}$. Since $\mathcal{H}^1(K)<+\infty$, there exists $\rho>0$ such that $\B_{x_0,C_0\rho}\subset\Om$ and $\mathrm{Card}(K\cap\partial \B_{x_0,\rho})<+\infty$. We denote by $K_0,K_1,\hdots,K_N$ the connected components of $K\cap \B_{x_0,\rho}$, with $K_0$ being the one that contains $x_0$. By compactness there exists $R\in (0,\rho)$ such that $K\cap\B_{x_0,R}\subset K_0$: thus for any $\xi,\zeta\in K\cap\B_{x_0,R}$, there exists some connected arc $\Gamma$ in $\overline{\B_{x_0,\rho}}$ joining $\xi,\zeta$ in $K_0$, and the chord-arc property applies (meaning $\mathcal{H}^1(\Gamma)\leq C_0|\zeta-\xi|$). Let now $x\in \B_{x_0,R}$, $r>0$ such that $\B_{x,r}\subset\B_{x_0,R}$, and assume that $\beta_K(x,r)\leq \eta$.

Without loss of generality we assume $x=0,r=1$, and assume that $\beta_K(x,r)=\beta_K(0,1)$ is reached for the affine line $\{z\in\R^2:z_2=0\}$, in the sense that
\begin{equation}\label{eq_twoside_hausdorff}
K\cap\B_{1}\subset\{z\in\R^2:|z_2|\leq\eta\},\ \{z\in\B_{1}:z_2=0\}\subset K+\overline{\B_\eta}.
\end{equation}
Taking $\left(-\frac{3}{4},0\right)$ which belongs to the left-hand side, there is an element $\zeta\in K$ such that $\left|\zeta-\left(-\frac{3}{4},0\right)\right|\leq \eta$. Let $K'$ the connected component of $\B_{1}\cap K(=\B_1\cap K_0)$ that contains $\zeta$. Let
$$s=\sup\{z_1,z\in K'\}$$
We have necessarily $s\geq -\frac{3}{4}-\eta$, and our goal is to prove that $s\geq \frac{1}{2}$: indeed since $K'\subset\{z:|z_2|\leq \eta\}$ this implies that $K'$ separates $D^\pm_{\frac{1}{2}}$ in $\B_{\frac{1}{2}}$.

Assume it is not the case, thus $s<\frac{1}{2}$. There exists $\zeta'\in K'$ such that $\zeta'_1=s$, and $|\zeta'_2|\leq\eta$. By \eqref{eq_twoside_hausdorff}, there exists $\zeta''\in \B_1\cap K\setminus K'$ such that $\left|\zeta''-\left(s+2\eta,0\right)\right|\leq \eta$, and we have
$$|\zeta'-\zeta''|\leq (1+\sqrt{5})\eta$$
However, since  $\B_1\cap K$ is a subset of $K_0$, there exists a connected arc $\Gamma\subset K_0\subset K\cap\B_{x_0, R}$ that connects $\zeta'$ to $\zeta''$, and since $\B_{x_0,C_0R}\subset\Om$, by the local chord-arc condition we must have
$$\mathcal{H}^1(\Gamma)\leq (1+\sqrt{5})C_0\eta.$$ 
Since $\Gamma$ is connected to $\zeta'$, which satisfies $|\zeta'_1|\leq \frac{3}{4}+\eta$, $|\zeta'_2|\leq \eta$, thus it is confined in the following ball:
$$\Gamma\subset\B_{\frac{3}{4}+(2+(1+\sqrt{5})C_0)\eta}.$$ 
With our definition of $\eta$ this implies that $\Gamma\subset\B_1$, which contradicts the fact that $\zeta',\zeta''$ belong to different connected components of $K\cap\B_1$. This ends the proof.
\end{proof}

\subsection{Epsilon-regularity}

The purpose of this section is to prove the following $\varepsilon$-regularity result.

\begin{proposition}\label{epsilonR} There exist some constants $\varepsilon>0$ depending only on $\A$ such that the following holds. Let $(u,K)$ be a topological Griffith minimizer in $\Omega \subset \R^2$ with $K$ connected, and let $x_0\in K$. Then there exists $R>0$ (depending on $K,\Om,x_0$) such that for all $x\in K\cap \B_{x_0,R}$ and $r>0$ such that $\B_{x,r}\subset \B_{x_0,R}$, if
$$ \beta(x,r) \leq \varepsilon,$$
then $K \cap \B_{x,\frac{1}{2}r}$ is a $C^{1,\alpha}$ regular curve.\end{proposition}
More precisely, we prove that if $\beta(x,r)(<\eps)$ is reached for a horizontal line, then $K\cap\B_{x,\frac{1}{2}r}$ is a graph in the first coordinate.

A first ingredient of the proof only uses the elasticity equation satisfied by $u$, with the polynomial decrease rate for the energy of the Airy function proved in \ref{prop_poly_growth}.
\begin{proposition}\label{prop_polygrowth_griffith}
There exists $C_1,\kappa>0$ depending only on $\A$ such that the following holds. Let $(u,K)$ be a topological Griffith minimizer in $\Om(\supset\B_1)$ such that $K$ is connected, $0\in K$, $$\beta(0,1)\leq \frac{1}{2}$$
and $K$ separates $D_{0,1}^{\pm}$ in $\B_{0,1}$. Then, for any $r\in \left[2\beta(0,1),1\right]$, we have
$$\om(x,r)\leq C_1 r^{\kappa}\om(x,1).$$
\end{proposition}
\begin{proof}
By the separation condition, for any such $r$, no connected component of
$$\left\{\theta\in\frac{\R}{2\pi\Z}: re^{i\theta}\notin K\right\}$$
has length larger than $\frac{4\pi}{3}$. Thus by Proposition \ref{prop_poly_growth}, denoting $\phi$ the Airy function of $u$, we have
$$\frac{1}{r}\int_{\B_r}|\nabla^2\phi|^2\leq C_{\frac{4\pi}{3}}r^{\kappa_{\frac{4\pi}{3}}}\int_{\B_1}|\nabla^2\phi|^2$$
Since $|\nabla^2\phi|^2=|\A e(u)|^2$, there exists some constant $C_\A>0$ such that
$$C_\A^{-1} |e(u)|^2\leq |\nabla^2\phi|^2\leq C_\A |e(u)|^2$$
and we get the result with $\kappa=\kappa_{\frac{4\pi}{3}}$, $C_1=C_\A^{2}C_{\frac{4\pi}{3}}$.
\end{proof}

The next ingredient in the proof of Proposition \ref{epsilonR} is the flatness estimate that was proven in \cite[prop. 3.2]{BIL} and that we recall here below (notice that we do not need to assume $K$ connected in the following statement). Moreover, it may be verified that the proof of this statement only relies on competitors which separate $D^\pm_{x,r}$ in $\B_{x,r}$, and thus applies to topological minimzers.

\begin{proposition}[\cite{BIL}, Proposition 3.2] \label{prop3.2.BIL}
There exist $\eps_2>0$ and $C_2>0$ (only depending on $\theta_0$, the Ahlfors regularity constant of $K$) such that the following property holds. Let $(u,K)$ be a minimizer of the Griffith functional in $\Om(\subset\R^2)$. For all $x_0 \in K$ and $r>0$ such that $\overline \B_{x_0,r}\subset \Omega$,
$$\omega(x_0,r)+\beta(x_0,r) \leq \eps_2,$$ 
and $K$ separates $D^\pm_{x_0,r}$ in $\overline{\B_{x_0,r}}$, we have
$$\beta\left(x_0,\frac{r}{50}\right) \leq  C_2 \omega(x_0,r)^{\frac{1}{14}}.$$
\end{proposition}

For the recall,  this statement was used in \cite{BIL}, together with an energy decay result (which was proved by compactness in \cite{BIL}), in order to obtain a $C^1$ regularity result on a connected component of $K$. Here we obtained the energy decay from the monotonicity by Proposition \ref{prop_polygrowth_griffith} instead of the compactness argument of \cite{BIL}. This leads to the following alternative proof of $\varepsilon$-regularity result for sufficiently flat and separating Griffith minimizers.
\begin{proposition}[$\varepsilon$-regularity]\label{prop_eps_reg_separating}
There exists $\eps>0$ depending only on $\A$ such that the following holds. Let $(u,K)$ be a topological Griffith minimizer in $\Om(\subset\R^2)$ such that $K$ is connected. Let $x\in K$, $r>0$ such that $\B_{x,r}\subset\Om$. Assume that
$$\beta(x,r)\leq \eps$$
and that $K$ separates $D^+_{x,r}$ and $D^-_{x,r}$. Then $K\cap \B_{x,\frac{1}{2}r}$ is a $\mathcal{C}^{1,\alpha}$ regular curve.
\end{proposition}

\begin{proof}

Assume without loss of generality that $r=1$. Since $K$ is connected, then $u$ admits a Airy function $\phi$ that belongs to $H^2_{0,K,\loc}(\Om,\R)$.

Consider $C_A$ the constant from Proposition \ref{prop_Ahlfors}, $C_1,\kappa>0$ from Proposition \ref{prop_polygrowth_griffith}, $C_2,\eps_2>0$ from Proposition \ref{prop3.2.BIL}.

We let
$$\tau=\min\left(\frac{1}{50},\left(2C_1^\frac{1}{14}\left(1+50^\frac{\kappa}{14}C_2\right)\right)^{-\frac{14}{\kappa}}\right),\  \eps_3=\min\left(50\tau\eps_2,\frac{1}{2}\tau\right),\ \overline{r}=\left(\frac{\eps_3}{2C_1C_A}\right)^{14\kappa^{-1}},\ \eps=\frac{1}{2}\overline{r}\eps_3$$
and we assume $\beta(x,1)\leq \eps$.

First, note by Proposition \ref{prop_polygrowth_griffith}, by the fact that $\om(x,1)\leq C_A$, and by the scaling properties of $\beta$, that for all $r\in (2\beta(x,1),1)$,we have
$$\om(x,r)\leq C_1 C_A r^\kappa,\ \beta(x,r)\leq r^{-1}\beta(x,1)$$
and that $K$ separates $D^\pm_{x,r}$ in $\B_{x,r}$. By our assumption on $\beta(x,1)$, we have $\overline{r}\geq 2\beta(x,1)$, and By definition of $\overline{r}$, this implies
$$\om(x,\overline{r})^\frac{1}{14}+\beta(x,\overline{r})\leq \eps_3$$
and $K$ separates $D^\pm_{x,\overline{r}}$ in $\B_{x,\overline{r}}$. Now, we prove by induction that for any $k\in\N$, we have
\begin{equation}\label{eq_induction_flatness}
\om(x,\tau^k\overline{r})^\frac{1}{14}+\beta(x,\tau^k\overline{r})\leq 2^{-k}\eps_3
\end{equation}
and $K$ separates $D^\pm_{x,\tau^k\overline{r}}$ in $\B_{x,\tau^k\overline{r}}$.
For $k=0$ this is already proved. Assume now that this is verified for some $k\in\N$. We have
$$\om(x,50\tau^{k+1}\overline{r})+\beta(x,50\tau^{k+1}\overline{r})\leq (50\tau)^{-1}2^{-k}\eps_3\leq \eps_2$$
and $\beta(x,\tau^k \overline{r})\leq \eps_3\leq \frac{1}{2}\tau$, so $D^{\pm}_{x,\tau^{k+1}\overline{r}}$ separates $\B_{x,\tau^{k+1}\overline{r}}$ and using propositions \ref{prop3.2.BIL} and \ref{prop_polygrowth_griffith} we get
\begin{align*}
\om(x,\tau^{k+1} \overline{r})^\frac{1}{14}+\beta(x,\tau^{k+1} \overline{r})&\leq C_1^\frac{1}{14}\tau^{\frac{\kappa}{14}}\om(x,\tau^k\overline{r})^\frac{1}{14}+C_2\om(x,50\tau^{k+1} \overline{r})^\frac{1}{14}\\
&\leq C_1^\frac{1}{14}\tau^{\frac{\kappa}{14}}\om(x,\tau^k\overline{r})^\frac{1}{14}+C_1^\frac{1}{14}C_2(50\tau)^{\frac{\kappa}{14}}\om(x,\tau^{k}\overline{r})^\frac{1}{14}\\
&\leq 2^{-(k+1)}\eps_3\text{ by def. of }\tau
\end{align*}

Thus, if $\beta(x,1)$ is sufficiently small, then for every $r\in (0,1)$, we have
$\beta(x,r)\leq C r^\alpha$ for some universal $C,\alpha>0$, and $K$ separates $D^{\pm}_{x,r}$ in $\B_{x,r}$.

By the standard regularity criterion (see for instance~\cite[Lemma~6.4]{BIL}), we conclude that if $\beta(x,1)$ is sufficiently small, then $K\cap \B_{x,\frac{1}{2}}$ is a $C^{1,\alpha}$ regular curve. This concludes the proof.
\end{proof}

We finally obtain Proposition \ref{epsilonR}.

\begin{proof}[Proof of Proposition \ref{epsilonR}] $R$ is given by Proposition \ref{propositionSepare}, and the result follows from Proposition \ref{prop_eps_reg_separating}. 
\end{proof}

\section{Proof of the main result}\label{sec_proof_main}

The purpose of this section is to prove Theorem \ref{main1} and \ref{main2} announced in the introduction. Our strategy is as follows: relying on the $\eps$-regularity and the monotonicity formula, we first prove Theorem \ref{main2} by a blow-up analysis. Theorem \ref{main1} is then obtained, by combining the regularity obtained in Theorem \ref{main2}, and an analysis of blow-down. 

\begin{proof}[Proof of Theorem \ref{main2}]
We denote by $\phi$ the Airy function associated to $u$, such that $\phi\in H^2_{0,K}(\Om,\R)$. Let $x\in K$, we denote 
$$\mathbf{M}_0:=\lim_{r\to 0}\mathbf{M}(\phi(x+\cdot),r).$$
We consider the full blow-up sequence  $(u_{x,r},K_{x,r})$ for $r\to 0$. By Proposition \ref{prop_compactness}, there exists some subsequence $r_n\to 0$ such that  $(u_{x,r_n},K_{x,r_n})$ converges to some  limit $(v,L)$. The Airy function of $u_{x,r_n}$ is $\phi_{x,r_n}:=r_n^{-\frac{3}{2}}\phi(x+r_n\cdot)$. $\phi_{x,r_n}$ is locally uniformly bounded in $H^2_\loc(\R^2,\R)$ by the upper Ahlfors regularity (see Proposition \ref{prop_Ahlfors}), and the continuity at the limit of the elastic energy from Proposition \ref{prop_compactness}, implies the strong $H^2_\loc(\R^2,\R)$ convergence of $\phi_{x,r_n}$ to the Airy function of $v$, denoted $\psi$ (see also \cite[Prop. 6.1]{BCL1} for a more detailed proof of convergence).
As a consequence, for any $r>0$, we have
$$\mathbf{M}(\psi,r)=\lim_{n\to \infty}\mathbf{M}(\phi(x+\cdot),r_n r)=\mathbf{M}_0.$$
By Theorem \ref{th_monotonie_disk}, this leaves two options: either $\mathbf{M}_0=0$, $\psi=0$ (meaning that $v$ is locally a rigid motion in $\R^2\setminus L$), and $L$ is a topological minimizer of its length, so it is either a line or a triple junction (not necessarily centered at the origin), or $\mathbf{M}_0>0$, $\psi$ is of the form given in \eqref{eq_description32homogeneous} for some $a>0$, and $L$ is a half-line starting from the origin.

We differentiate the three cases:\bigbreak

\textbf{$v$ is locally rigid and $L$ is a line}. In particular, $\beta_K(x,r_n)\to0$. For $n$ large enough, Proposition~\ref{epsilonR} applies and shows that $K$ is a $C^{1,\alpha}$ curve in a neighborhood of $x$. It follows in particular that every blow-up at $x$ is the tangent line to this curve. Moreover, after possibly reducing the neighborhood, the same conclusion holds at every point of $K$ in this neighborhood. Hence $\mathcal R$ is relatively open in $K$, and neither a triple junction nor cracktip can occur as another blow-up at a point of $\mathcal R$.\bigbreak

\textbf{$v$ is locally rigid and $L$ is a triple junction}. If the triple junction is not centered, then an additional blow-up proves that for some $\eps_n\to 0$, the blow-up sequence $(u_{x,\eps_n r_n},K_{x,\eps_n r_n})$ converges to a rigid motion and a line: this is excluded by the previous point.

Thus for any sequence $s_n\to 0$, there is some  subsequence such that the associated blow-up sequence converges to a (possibly different) triple junction. As a consequence, there exists some $r_{x}>0$ such that the following holds: for any $r\in (0,2r_x)$, the set
$$
K_{x,r}\cap \B_{0,3/4}\setminus \overline{\B_{0,1/4}}
$$
is contained in three disjoint thin neighborhoods of the three arms of a $Y$-cone.

Let $z$ be a point on one of these arms, away from the origin. In a small ball centered at $z$, the $Y$-cone coincides with a straight line.
Consequently, if $y\in K$ corresponds after rescaling to a point close to $z$, then, for a fixed sufficiently small constant $c>0$, $\beta_K(y,cr)$ is arbitrarily small when $r$ is small. Proposition~\ref{epsilonR} then implies that $K\cap \B_{y,acr}$ is a $C^{1,\alpha}$ regular curve.

Covering each of the three arms in the annulus by finitely many such balls, we obtain that, for every sufficiently small $r$,
\[
    K\cap
    \left(
        \B_{x,3r/4}\setminus \overline{\B_{x,r/4}}
    \right)
\]
is the disjoint union of three $C^{1,\alpha}$ arcs, one near each arm of the corresponding $Y$-cone. It is important here that Proposition~\ref{epsilonR} applies to the whole set $K$ in the smaller ball. Thus no additional piece of $K$ can be hidden inside one of the three tubular neighborhoods.

Applying the same argument at successive smaller scales and using the overlap of the corresponding annuli, the three arcs glue together. Indeed, in every overlap $K$ is locally a single $C^{1,\alpha}$ curve, so a branch cannot stop, split, or merge with another one. We conclude that
\[
    (K\cap \B_{x,r_x})\setminus\{x\}
\]
has exactly three connected components, each of which is a $C^{1,\alpha}$ arc away from $x$ and has $x$ in its closure. For any $r\in (0,r_x]$, we may define $\theta_1(r)<\theta_2(r)<\theta_3(r)(<\theta_1(r)+2\pi)\in\R$, that has a $\mathcal{C}^{1,\alpha}$ dependency in $r$, such that $x+re^{i\theta_k(r)}$ belongs to the $k$-th arm (for some arbitrary choice). An homeomorphism $\B_{x,r_x}\to \B_{r_x}$ as in the statement of Theorem \ref{main2} is then given by
$$\Phi(x+re^{i\theta})=\begin{cases}
re^{\frac{2\pi}{3}i\frac{\theta-\theta_1(r)}{\theta_2(r)-\theta_1(r)}}&\text{ if }\theta_1(r)\leq \theta\leq \theta_2(r)\text{ mod }2\pi\\
re^{\frac{2\pi}{3}i\frac{\theta-\theta_2(r)}{\theta_3(r)-\theta_2(r)}+\frac{2\pi}{3}i}&\text{ if }\theta_2(r)\leq \theta\leq \theta_3(r)\text{ mod }2\pi\\
re^{\frac{2\pi}{3}i\frac{\theta-\theta_3(r)}{\theta_1(r)+2\pi-\theta_3(r)}+\frac{4\pi}{3}i}&\text{ if }\theta_3(r)\leq \theta\leq \theta_1(r)+2\pi\text{ mod }2\pi\\
\end{cases}.$$
for $r\in (0,r_x]$.\bigbreak

\textbf{$\mathbf{M}_0>0$, $\psi$ is of the form given in \eqref{eq_description32homogeneous} for some $a>0$, and $L$ is a half-line starting from the origin.} Up to a rotation, this means that $(v,L)=(a+k_{sl}u_{sl}+k_{op}u_{op},\R_-\times\{0\})$, for some rigid motion $a$, and some coefficients $(k_{op},k_{sl})\in\R^2\setminus\{(0,0)\}$. By Proposition \ref{prop_instability_slide}, we have exactly $(k_{op},k_{sl})=(1,0)$ or $(-1,0)$. The rest follows the same argument as the previous case: for some sufficiently small $r_{x}>0$, we have that for all $r\in (0,2r_x)$ the set
$$\frac{K-x}{r}\cap\B_{0,\frac{3}{4}}\setminus\B_{0,\frac{1}{4}}$$
is contained in a thin neighbourhood of a radial segment connecting $\partial\B_{0,\frac{1}{4}}$ and $\partial\B_{0,\frac{3}{4}}$, so that by application of the $\eps$-regularity result on a covering of this segment, we find that $K\cap\B_{x,r_x}\setminus\{x\}$ is a $\mathcal{C}^{1,\alpha}$ curve that meets once every concentric circle and ends at $x$. Denoting a smoothly-defined $\theta(r)\in\R$ such that $x+r e^{i\theta(r)}\in K$, the homeomorphism of Theorem \ref{main2} is given by
$$\Phi(x+re^{i\theta})=re^{i(\theta-\theta(r)+\pi)}.$$
\end{proof}

To prove the classification of global minimizers, we rely on a similar blow-\textbf{down} argument. We begin by recalling the following result, established in \cite{LLL}, which asserts that $\mathbb{R}^2\setminus K$ has no bounded connected components. This property is quite natural ; suppose that $D\subset \mathbb{R}^2\setminus K$ is a bounded connected component. By minimality, $u$ must coincide with a rigid motion in $D$. One can then construct a competitor by suitably ``opening'' the component $D$, thereby producing a sufficiently large decrease in the energy to contradict the minimality of $(u,K)$. This is essentially the argument used in the proof of \cite{LLL}.

\begin{proposition}[No isolated component]\label{prop_boundedcomponent}{\rm \cite[Proposition 3.1]{LLL}} Let $(u, K)$ be a global minimizer of the Griffith functional. Then $\R^2\setminus K$  has no bounded connected component.
\end{proposition}

We are now ready to prove the classification of global minimizers with connected crack.

\begin{proof}[Proof of Theorem \ref{main1}] Let $(u,K)$ be a global minimizer in $\R^2$ with $K$ connected. Let us first notice that $K$ is not bounded. Otherwise if $K\subset \B_{x_0,R_0}$ with $x_0\in K$, then by  lower Ahlfors regularity of $K$ (see \cite[Th. 2.3]{FLK2}) we would have, for all $R>R_0$, for some universal constant $C>0$,
$$CR\leq \mathcal{H}^1(K\cap \B_{x_0,R})=\mathcal{H}^1(K\cap \B_{x_0,R_0})\leq C_A R_0,$$
a contradiction when $R\to +\infty$. In particular we deduce that for all $x\in K$ and $R>0$, $K\cap \partial\B_{x,R}\neq \emptyset$.

We let  $\phi$ be the Airy function associated to $u$. We let
$$\mathbf{M}_\infty=\lim_{r\to +\infty}\mathbf{M}(\phi,r)$$
Note that this quantity is invariant by coordinate translation of $\psi$, since $\psi_{x,r}-\psi_{y,r}$ converges in $H^2_\loc(\R^2,\R)$ to $0$ as $r\to +\infty$, for any fixed, $x,y$.

If $\mathbf{M}_\infty=0$, then by Theorem \ref{th_monotonie_disk}, for any $x\in K$, for any $r>0$, we have $0\leq \mathbf{M}(\phi(x+\cdot),r)\leq \mathbf{M}_\infty=0$, so $\phi\equiv 0$. As a consequence, $u$ is locally rigid, so $K$ is either a line or a triple junction.

Assume now that $\mathbf{M}_\infty>0$. Let  $x\in K$. Following the same blow-up argument as the proof Theorem \ref{main2}, there is some sequence $r_n\to +\infty$ such that
$$(u_{x,r_n},K_{x,r_n})\to (v,L)$$
and $\phi_{x,r_n}\to \psi$ in $H^2_\loc(\R^2,\R)$, where $(v,L)$ is a global minimizer, and $\psi$ is the Airy function of $v$. We have for every $r$ that
$$\mathbf{M}(\psi,r)=\lim_{n\to +\infty}\mathbf{M}(\phi(x+\cdot),r_nr)=\mathbf{M}_\infty(>0),$$
so $\psi$ is $\frac{3}{2}$-homogeneous and non-zero: arguing as earlier, by Proposition \ref{prop_instability_slide}, $(v,L)$ is (up to a rotation of the coordinates) of the form $(a\pm u_{op},\R_-\times\{0\})$ for some rigid motion $a$, so in particular
$$\mathbf{M}_\infty=\mathbf{M}(\psi,1)=\mathbf{M}(\phi_{op},1).$$

For some sufficiently large $R_x>0$, we then have that for all $R\geq R_x$, the set
$$\frac{K-x}{R}\cap\B_{0,\frac{3}{4}}\setminus\B_{0,\frac{1}{4}}$$
is contained in a thin neighbourhood of a radial segment joining $\partial\B_{0,\frac{3}{4}}$ and $\partial\B_{0,\frac{1}{4}}$, so that by application of the $\eps$-regularity result \ref{prop_eps_reg_separating} on a covering of this segment, we obtain (following the same covering argument as the proof of Theorem \ref{main2}) that for some $\mathcal{C}^1$ function $\theta:[R_x,+\infty)\to\R$, we have
$$K\setminus\B_{x,R_x}=\{x+Re^{i\theta(R)}, R\geq R_x\}.$$
By Theorem \ref{main2}, the sets $\mathcal{Y}$ and $\mathcal{C}$ are finite, and we consider the planar graph $G$ formed by the vertices $\mathcal{Y}\cup\mathcal{C}\cup\{\infty\}$ with edges $\mathcal{R}$. The vertices of $\mathcal{Y}$ have degree $3$ and every other vertex has degree $1$. 

Since $\R^2\setminus K$ has no bounded component by Proposition \ref{prop_boundedcomponent}, $G$ has no cycle, so it is a tree, meaning its set of leaves $\mathcal{C}\cup\{\infty\}$ must have at least two elements. Thus $\mathcal{C}$ is non-empty and we denote $z\in\mathcal{C}$ some cracktip point.

Since any blow-up of $\phi$ at $z$ is some rotation of $\phi_{op}$ we then have, for any $r>0$,
$$\mathbf{M}(\phi_{op},1)=\lim_{s\to 0}\mathbf{M}(\phi(z+\cdot),s)\leq \mathbf{M}(\phi(z+\cdot),r)\leq \lim_{s\to +\infty}\mathbf{M}(\phi(z+\cdot),s)=\mathbf{M}(\phi_{op},1).$$
Thus this is a chain of equalities and $\phi(z+\cdot)$ is a nonzero, $\frac{3}{2}$-homogeneous biharmonic function: by Proposition \ref{prop_instability_slide}, we obtain that $(u,K)$ is of the  form \eqref{eq_def_cracktip}.
\end{proof}

\vspace{1cm}
\noindent {\bf Acknowledgment.}  The authors' research was supported by the project ANR STOIQUES (ANR-24-CE40-2216) financed by the French Agence Nationale de la Recherche (ANR).

\bibliographystyle{plain}
\bibliography{biblio}

\end{document}